\documentclass[a4paper,10pt]{article}
\usepackage{tikz}
\usepackage[utf8]{inputenc}
\usepackage{amsthm}
 
\newtheorem{proposition}{Proposition}[section]
\newcommand{\mycomment}[1]{ }

\def\Oij{\bO_{i,j}}

\usepackage{amssymb}
\usepackage{amsfonts}
\usepackage{amsmath}

\newcommand\beq{\begin{equation}}
\newcommand\eeq{\end{equation}}

\usepackage{dsfont}

\newcommand{\dt}{\Delta t}

\renewcommand{\to}{\rightarrow}

\renewcommand{\det}{\operatorname{det}}
\renewcommand{\div}{\operatorname{div}}

\newcommand{\tr}{\operatorname{tr}}

\newcommand{\spn}{\operatorname{Span}}

\newcommand{\grad}{\boldsymbol{\nabla}}
\renewcommand{\div}{\operatorname{div}}

\newcommand{\bB}{\boldsymbol{B}} 

\newcommand{\bC}{\boldsymbol{C}}

\def\bD{\boldsymbol{D}}

\newcommand{\be}{\boldsymbol{e}}

\newcommand{\bF}{\boldsymbol{F}}
\newcommand{\bg}{\boldsymbol{g}}
\newcommand{\bG}{\boldsymbol{G}}

\newcommand{\bI}{\boldsymbol{I}} 

\newcommand{\bL}{\boldsymbol{L}}

\def\bn{\boldsymbol n}

\newcommand{\bO}{\boldsymbol{O}}

\newcommand{\bP}{\boldsymbol{P}}

\def\bR{\boldsymbol R}

\newcommand{\bu}{\boldsymbol{u}}
\newcommand{\bU}{\boldsymbol{U}}
\newcommand{\gbu}{{\grad\bu}}

\newcommand{\bW}{\boldsymbol{W}}

\newcommand{\bOmega}{\boldsymbol{\Omega}} 

\newcommand{\bSigma}{\boldsymbol{\Sigma}}

\newcommand{\Strs}{{\boldsymbol{\Sigma}}}

\def\Ecal{\mathcal{E}} 			

\newcommand{\NN}{\mathbb{N}}

 \newcommand{\R}{\mathbb{R}} 

\title{
Derivation and numerical approximation of
hyperbolic viscoelastic flow systems:
Saint-Venant 2D equations for Maxwell fluids 
}
\author{
S\'ebastien Boyaval, Laboratoire d'hydraulique Saint-Venant \\
(EDF R\&D -- Ecole des Ponts ParisTech -- CEREMA) Universit\'e Paris-Est \\
EDF'lab 6 quai Watier 78401 Chatou Cedex France, \\
\& INRIA Paris, MATHERIALS (sebastien.boyaval@enpc.fr)
}
\newcommand{\ro}{h}
\newcommand{\ux}{u}
\newcommand{\uy}{v}

\newcommand{\sx}{c_{xx}}
\newcommand{\sy}{c_{yy}}
\newcommand{\sxy}{c_{xy}}
\newcommand{\sz}{c_{zz}}

\begin{document}
\maketitle
\begin{abstract}

We pursue here the development of models for complex (viscoelastic) fluids in \emph{shallow free-surface gravity flows}
which was initiated by [Bouchut-Boyaval, M3AS (23) 2013] for 1D (translation invariant) cases.

\smallskip

The models we propose are \emph{hyperbolic} quasilinear 
systems that generalize Saint-Venant shallow-water equations to incompressible 
\emph{Maxwell fluids}.
The models are compatible with 
a formulation of the thermodynamics second principle.
In comparison with Saint-Venant standard shallow-water model, the momentum balance 
includes extra-stresses 
associated with an elastic potential energy in addition to a hydrostatic pressure.
The extra-stresses are determined by an additional tensor 
variable solution to a differential 
equation with various possible time rates. 

\smallskip

For the numerical evaluation of 
solutions to Cauchy problems, 
we also propose \emph{explicit} schemes discretizing our generalized Saint-Venant systems  
with Finite-Volume approximations that are \emph{entropy-consistent} (under a CFL constraint)
in addition to satisfy exact (discrete) mass and momentum conservation laws.
In comparison with most standard viscoelastic numerical models,
our discrete models can be used for any retardation-time values
(i.e. in the vanishing ``solvent-viscosity'' limit).

\smallskip

We finally illustrate our hyperbolic viscoelastic flow models numerically
using computer simulations in benchmark test cases. 
On extending to Maxwell fluids some free-shear flow testcases that are standard benchmarks for Newtonian fluids,
we first show that our (numerical) models reproduce well the viscoelastic physics, phenomenologically at least,
with zero retardation-time.
Moreover, with a view to 
quantitative 
evaluations, numerical results in the lid-driven cavity testcase show that, in fact, 
our models can be compared with standard viscoelastic flow models 
in sheared-flow benchmarks on adequately choosing the physical parameters of our models.
%
Analyzing our models 
asymptotics should therefore shed new light on the famous
High-Weissenberg Number Problem (HWNP), 
which is a limit for all the existing viscoelastic numerical models. 

\mycomment{
To cope with the limitations of Maxwell \emph{linear} and \emph{infinitely-extensible} viscoelastic rheology [] 
we propose here to study finitely-extensible nonlinear viscoelastic fluids.
}
\end{abstract}

\section{Introduction} 
\label{sec:intro}

Modelling the viscoelastic large deformation 
of a flowing continuum material is still a challenge, despite more than a century of attempts after Maxwell \cite{Maxwell01011867}. 
%
%
Precisely, to numerically predict the evolution in time of realistic non-Newtonian \emph{fluid} continuum materials, 
one does not know yet a good mathematical model for 
Eulerian fields 
characterizing a multidimensional flow 
that would properly generalize Maxwell's one-dimensional (1D) model. 

In addition to respect the fundamental physical principles (Galilean invariance, mass conservation etc.), a good viscoelastic model should account
i) for the stresses generated within the fluid by a loss of elastic (potential) energy in deformations -- shear and extension --, 
and ii) for the (viscous) energy dissipation leading the fluid to a 
flow equilibirium.

Many \emph{constitutive equations} have been proposed 
to relate the stresses with the symmetrized velocity gradient as a measure of strain,
which are formally satisfying from the loose viewpoint above,
but research is still going on \cite{bird-drugan-2017}.
Indeed, 
there seems to be no general model yet that allows one to correctly match the 
experimentally-observed 
behaviour of complex fluids in steady flows, like shear-thinning and strain-hardenning,
see e.g. \cite{bird-curtiss-armstrong-hassager-1987a,renardy-2000}.

Moreover, in any case, a good model should also allow one to predict the evolution in time of a piece of fluid, 
at least for small times, knowing initial conditions (a Cauchy initial-value problem).
Now, 
most numerical viscoelastic flow model 
use similar sets of 
PDEs 
for a mass density $\rho$, a \emph{solenoidal} velocity 
$\bu$, and a 
2-tensor 
$\Strs$ with various possible physical interpretations with the Eulerian viewpoint. 
Typically, differential models of \emph{rate-type} have been proposed for Maxwell fluids 
when the tensor field $\Strs$ is either Cauchy's stress tensor 
or a strain tensor like Cauchy-Green, 
see e.g. \cite{oldroyd-1950,bird-curtiss-armstrong-hassager-1987a,renardy-2000}.
But despite some recent mathematical progress showing that Cauchy problems 
for 
physically-interesting models like FENE-P and Giesekus 
have global solutions \cite{masmoudi-2011}, 
the \emph{numerical computation} of solutions remains unsatisfying too.


There has been continuing progress in the numerical discretization of Cauchy problems for viscoelastic models
\cite{owens-philips-2002,fattal-kupferman-2005,PAN2007283,afonso-oliveira-pinho-alves-2009,knetchges-behr-elberti-2014,Sousa2016129}
and its analysis \cite{boyaval-lelievre-mangoubi-2009,barrett-boyaval-2010,perrotti-walkington-wang-2017,barrett-boyaval-2017}.
However, numerical computations remain impossible in some physically-meaningful configurations 
at moderately high Weissenberg number in particular, 
which is quite frustrating for applications.

We note that 
all the discrete 
viscoelastic models 
mentioned above
consider the \emph{solenoidal} flows of incompressible Maxwell fluids, 
with a non-zero \emph{retardation time}
(i.e. they have a ``background'' viscous stress term justified by the \emph{solvent viscosity} in the context of polymer flows).
Then, they rely on the assumption that velocity gradients remain bounded,
and numerical \emph{steady-state flows} are usually computed
using iterative algorithms to approximate the nonlinear terms (of the constitutive equations at least). 
Now, numerical instabilities typically occur when the nonlinear terms prevail in the (discrete) equations
and iterative algorithms fail to converge (to an implicit time-discretization e.g.).
%

\mycomment{
Is it possible to numerically simulate standard viscoelastic flows without background viscosity, and reproduce instabilities or not ?
Note, in practice, we still need a viscous term to impose BC and reach the -- same -- large-time viscous equilibrium as in standard benchmarks.
But our free-energy estimate does not require that term, to hold e.g. under bounded velocity gradients !
}

In this work, we would like to study 
discrete 
viscoelastic models 
that are explicitly computable
(no iterative algorithm is required for the approximation of nonlinearities 
as fixed points)
and that do not assume a non-zero retardation time (velocity gradients are possibly undounded).
Hence, we consider quasilinear systems of first-order with possibly non-solenoidal velocities,
that express mass and energy conservation for 
elastic fluids with one 
single relaxation time (so-called Maxwell fluids). 
Moreover, for the sake of stability, we consider hyperbolic systems endowed with an additional conservation law
(in smooth evolutions) as a formulation of the second thermodynamics principle.
Even if one single additional conservation law may not suffice to define \emph{unique} admissible solutions
to multidimensional hyperbolic systems, see e.g. \cite{chiodaroli-delellis-kreml-2015},
it can be useful for numerical stability purposes at least,
see e.g. 
\cite{boyaval-lelievre-mangoubi-2009,barrett-boyaval-2010,perrotti-walkington-wang-2017,barrett-boyaval-2017}.

A quasilinear system of PDEs has been considered as a model for viscoelastic 2D flows of a slightly compressible Maxwell fluid 
in \cite{PHELAN1989197,EDWARDS1990411,Guaily2010158,Guaily2011258}.
It uses the Upper-Convected Maxwell equation for the stress tensor $\Strs$ (like Oldroyd-B model) 
and it is hyperbolic under a simple physically-natural condition 
\emph{that is independent of the system 
parameters} when rewritten 
in terms of a deformation tensor (the left Cauchy-Green deformation tensor, see below in Section \ref{sec:saintvenant}).
Numerical computations show the interest 
of the approach. 
However, the model 
does not conserve mass 
and
is not obviously compatible with thermodynamics principles.
Moreover, 
it does not seem easily generalized to 3D flows,
and its stress-strain relationship is physically unclear.

\mycomment{
We also refer to \cite{Liapidevskii2013} for a special class of solutions using hyperbolic \emph{submodels} of standard system of PDEs for viscoelastic flows.
}

\medskip

In this work, we use {the simplified Saint-Venant framework} for hydrostatic 2D (shallow) free-surface flows 
to propose and study 
quasilinear systems such that:
(i) they model incompressible viscoelastic flows of Maxwell fluids with a \emph{deformation tensor} as dependent variable,
(ii) they conserve mass and satisfy a formulation of the thermodynamics principles,
(iii) they remain hyperbolic under simple physically-natural admissibility conditions.
%
Two possible models are identified
in Section~\ref{sec:model}.

One model uses the 
Upper-Convected Maxwell equation for a deformation tensor of Cauchy-Green left type.
It is very similar to that in \cite{PHELAN1989197,EDWARDS1990411,Guaily2010158,Guaily2011258}
except that it is rewritten in terms of a deformation tensor
(which then allows us to show why 
it is hyperbolic under a 
condition independent of the system parameters).
In comparison, our version additionally conserves mass 
(to that aim, we introduce a free-surface and decompose the pressure in two terms) 
and it is endowed with a formulation of the thermodynamics principles.
This is one of the 2D 
models for viscoelastic incompressible shallow free-surface 
flows under gravity
that generalize Saint-Venant's approach (initially for shallow \emph{water}) to non-Newtonian 
fluids, which we previously derived 
in 
\cite{bouchut-boyaval-2015} after 
depth-averageing Oldroyd-B model 
for fast thin-layer flows with small 
viscosity.

Another model is proposed that uses \emph{Finger} deformation tensor (the \emph{inverse} of Cauchy-Green left one),
with a clear stress-strain relationship (using Euler-Almansi strain) 
and with \emph{Cotter-Rivlin} natural frame-invariant time-rate-of-change.
The model is hyperbolic under the same simple physically-natural conditions as the previous one,
it conserves mass,
and it is endowed with a similar formulation of the thermodynamics principles.
We also note that the model is similar to Reynolds-averaged models that have been proposed 
for the 
Reynolds stresses in weakly-sheared turbulent flows \cite{teshukov-2007,Teshukov2007}, see also \cite{hal-01527469}.

In Section~\ref{sec:FV}, we propose entropy-consistent Finite-Volume discretizations of the two models.
The numerical scheme is a 2D extension of the Suliciu-type relaxation approach developped in \cite{bouchut-boyaval-2013} 
for a closed subsystem fully describing the 1D (translation-invariant) motions 
(see also \cite{2016arXiv161108491B} for the exact solutions to the 1D Riemann problems).

In Section~\ref{sec:num}, we perform numerical simulations.
The numerical results show that our model are physically meaningful,
and a viable approach to the computer simulation of physically-realistic viscoelastic flows. 

\section{Hyperbolic 
viscoelastic Saint-Venant systems} 
\label{sec:model}

\subsection{Saint-Venant models for shallow free-surface flows}
\label{sec:saintvenant}

Given a constant gravity field $\bg=-g\be_z$, let us equip space with a Cartesian frame $(\be_x,\be_y,\be_z)$
and consider the flow, i.e. the evolution in time and space, of a viscoelastic fluid material
with a non-folded free-surface $z=H(t,x,y)\ge0$ above a flat impermeable plane $z=0$. 
We assume the fluid material incompressible with constant mass density.
Then, 
the dynamics of the flow is 
governed by a kinematic condition for the free-surface (in virtue of the mass conservation principle)
coupled to 
a 
momentum-balance equation for the 
velocity $(u,v,w)$
following 
physical principles.

Assuming 
the flow \emph{stratified} (i.e. with the 
acceleration small in vertical direction, 
and with the horizontal components $(u,v)$ of the velocity 
mostly uniform in vertical direction) such that
$$
H \int_0^H dz\, u^2 \approx \left(\int_0^H dz\, u\right)^2
\quad
H \int_0^H dz\, v^2 \approx \left(\int_0^H dz\, v\right)^2
$$
one standardly obtains an interesting system of quasilinear equations for $H>0$ 
and the depth-averaged velocity $\bU=(U,V)$ with components $U \approx \frac1H\int_0^H dz u$, $V \approx \frac1H\int_0^H dz v$
in the form of conservation laws:
\begin{align}
\label{eq:SVH}
& \partial_t H + \div( H \bU ) = 0
\\
\label{eq:SVHUV}
& \partial_t(H\bU) + \div(H\bU\otimes\bU + H(P+\Sigma_{zz})\bI - H\bSigma_h ) = -kH\bU
\end{align}
equivalently:
\begin{align*}
& \partial_t H + \partial_x( H U ) + \partial_y(H V) = 0
\\
& \partial_t(HU) + \partial_x( HUU + HP+H\Sigma_{zz}-H\Sigma_{xx} ) + \partial_y( HUV - H\Sigma_{xy} ) = -kHU
\\
& \partial_t(HV) + \partial_x( HUV - H\Sigma_{yx} ) + \partial_y (  HVV + HP+H\Sigma_{zz}-H\Sigma_{yy} ) = -kHV
\end{align*}
which is reminiscent of the 2D gas-dynamics equations of Euler 
with mass density proportional to $H$,
with \emph{specific} pressure $P+\Sigma_{zz}$
and with \emph{specific} 
stress-deviator $\Strs_h$. 
($P$ and $\Sigma_{ij}$ have the dimension of energy per unit mass 
rather than per unit volume, unlike standard Cauchy stresses and pressures typically measured in Pa,  
that is why we term them 
\emph{specific}. 
But for the sake of simplicity, we 
omit the label ``specific'' below, 
insofar as this is not ambiguous here.) 

With $P=gH/2$ the depth-averaged hydrostatic pressure, 
one retrieves the inviscid shallow-water model of Saint-Venant \cite{saint-venant-1871} when $\Sigma_{ij}=0$,
which coincides with Euler 
isentropic 2D flow model of perfect gases with $\gamma=2$.
Euler 2D model is (symmetric) hyperbolic 
and 
solutions to 
Cauchy problems allow one to numerically 
predict small-time evolutions 
of some compressible gases \cite{chen-wang-2002}. 
The shallow-water model of Saint-Venant is widely used e.g. in environmental hydraulics
for nonlinear \emph{irrotational} water waves in flat open channels with rugosity $k\ge 0$ \cite{chow-1959}, 
at least to predict the dynamics of \emph{long surface waves}.

Moreover, the system (\ref{eq:SVH}--\ref{eq:SVHUV}) can also sustain shear motions.
Formally, it 
approximates the depth-averaged free-surface Navier-Stokes equations 
for a viscous fluid with (small) 
viscosity $\nu\ge 0$ 
provided $\Sigma_{zz}=-(\Sigma_{xx}+\Sigma_{yy})$ 
and
\begin{equation}
\label{eq:viscoussteadystate}
\Sigma_{xx}=2\nu\partial_xU \quad \Sigma_{yy}=2\nu\partial_yV \quad \Sigma_{xy}=\nu(\partial_xV+\partial_yU)=\Sigma_{yx}
\end{equation}
see \cite{gerbeau-perthame-2001,marche-2007,bouchut-boyaval-2015} and references therein. 
%
%
%
An additional conservation law 
is 
satisfied:
\beq
\label{eq:energyconservation} 
\partial_t(HE) + \div(HE\bU+H(P+\Sigma_{zz})\bU-H\bSigma_h\cdot\bU) = -kH |\bU|^2 - HD \,,
\eeq
with viscous dissipation $D=2\nu\left(|D(\bU)|^2+2|\div\bu|^2\right)\ge0$ in shear,
for 
\begin{equation}
\label{mechanicalenergy} 
E = \frac{1}2 \left( |\bU|^2 + gH \right)
\end{equation}
the so-called \emph{mechanical energy}, convex in $(H^{-1},\bU)$. 
This is a formulation of the thermodynamics first principle 
when $k=0=\nu$ which implies, by 
Godunov-Mock theorem see e.g. \cite{godlewski-raviart-1996},
that the 
homogeneous system of 
conservation laws is (symmetric) \emph{hyperbolic}.
It also motivates the computation of small-time evolutions 
as $L^2$ solutions to 
Cauchy problems, 
see e.g.~\cite{serre1-1996}.
In particular, \emph{entropy solutions} to Cauchy problems 
exist that satisfy the inequality \eqref{eq:energyconservation}
\beq
\label{eq:energydissipation} 
\partial_t(HE) + \div(HE\bU+H(P+\Sigma_{zz})\bU-H\bSigma_h\cdot\bU) \le -kH |\bU|^2 - HD 
\eeq
as a formulation of the thermodynamics second principle. 
Moreover, the 1D (translation-invariant) bounded measurable entropy solutions are unique since the system is also strictly hyperbolic when $H>0$.

As soon as $\nu>0$ however (e.g. to account for shear motions -- 2D horizontal ones --), 
variations in the strain 
$D(\bU)=(\nabla\bU+\nabla\bU^T)/2$ have an \emph{infinite} propagation speed, 
and energy is dissipated everywhere in the strained fluid microstructure. 
This is not realistic for a number of flows, which may be better modelled as viscoelastic fluid continua of Maxwell type
(i.e. fluids endowed with an energy storage capacity of elastic modulus $G$ per unit mass 
and with a relaxation time $\lambda$ that define a viscosity $\nu=G\lambda$ together).

\subsection{Generalized Saint-Venant models for Maxwell fluids}
\label{sec:maxwell}

In addition to a viscosity $\nu$, viscoelastic fluids of Maxwell type 
are characterized by a time scale $\lambda$ 
for the 
relaxation of the stress to a viscous state proportional to strain
(it defines a Weissenberg number when compared with a time scale of the flow like $|\grad\bU|^{-1}$).
Moreover, an evolution equation is also needed for $\Sigma_{zz}$ and the horizontal stress \emph{tensor} 
$$
 \bSigma_h = \Sigma_{xx}\be_x\otimes\be_x + \Sigma_{yy} \be_y\otimes\be_y + \Sigma_{xy} \be_x\otimes\be_y + \Sigma_{yx} \be_y\otimes\be_x
 $$ 
to 
define a generalized Saint-Venant model for viscoelastic Maxwell fluids with $P=gH/2$. 
Various rate-type differential equations exist for the components of the 
tensor field $\bSigma_h$
that are compatible with Galilean invariance.

Following the same depth-averaged analysis as in \cite{gerbeau-perthame-2001,marche-2007} for viscous fluids, 
one can derive 
equations of the form 
$\lambda\stackrel{\triangle}{\bSigma} + \bSigma = 2\nu D(u,v,w)$ 
for the depth-averaged extra-stress variable $\bSigma=\bSigma_h+\Sigma_{zz}\be_z\otimes\be_z$ in a 3D Maxwell-fluid model,
using some 
admissible time-rate $\stackrel{\triangle}{\bSigma}$ 
(like the upper-convective derivative of Oldroyd-B model)
and a steady equilibrium $\bSigma \approx 2\nu D(u,v,w)$ 
(typically compatible with \eqref{eq:viscoussteadystate}). 
Depending on scaling 
assumptions, 
various equations 
arise in \cite{bouchut-boyaval-2015},
like (see \cite[section 6.1.3.]{bouchut-boyaval-2015}) 
\begin{equation}
\label{eq:Sigmah}
D_t \bSigma_h - \bL_h\bSigma_h - \bSigma_h \bL_h^T + \zeta (\bD_h\bSigma_h + \bSigma_h \bD_h^T) = (2\nu\bD_h-\bSigma_h)/\lambda
\end{equation}
where $D_t:=\partial_{t} + U\partial_{x} + V\partial_{y}\equiv \partial_{t} + \bU\cdot\grad$, 
$\bL_h:=\grad\bU$, $\bD_h:=(\bL_h+\bL_h^T)/2$, $\bW_h:=(\bL_h-\bL_h^T)/2$ for $\zeta=0$,
along with
\begin{equation}
\label{eq:Sigmazz}
D_t \Sigma_{zz} + 2(1-\zeta)\div\bU \Sigma_{zz} = (-2\nu\div\bU-\Sigma_{zz})/ \lambda.
\end{equation}
Equivalently, with $\Sigma_{xy}=\Sigma_{yx}$, (\ref{eq:Sigmah}--\ref{eq:Sigmazz}) reads
$$
\left\{
\begin{aligned}
D_t \Sigma_{xx} 
& 
- 
\left( 
2(1-\zeta)\Sigma_{xx} \partial_x U + \Sigma_{xy}\left((2-\zeta)\partial_yU-\zeta\partial_xV\right) 
\right) 
= (2\nu\partial_xU-\Sigma_{xx}) / \lambda 
\\
D_t \Sigma_{yy} 
&
- 
\left( 2(1-\zeta)\Sigma_{yy} \partial_y V + \Sigma_{xy}\left(-\zeta\partial_yU+(2-\zeta)\partial_xV\right) \right)
= (2\nu\partial_yV-\Sigma_{yy}) / \lambda 
\\
D_t \Sigma_{xy} 
&
- 
\Big( (1-\zeta/2)(\Sigma_{xx}\partial_xV+\Sigma_{yy}\partial_yU) + (-\zeta/2)(\Sigma_{yy}\partial_xV+\Sigma_{xx}\partial_yU)
\nonumber
\\
&
+ (1-\zeta)\Sigma_{xy}(\partial_xU+\partial_yV) \Big)
= (\nu(\partial_xV+\partial_yU)-\Sigma_{xy}) / \lambda 
\\
\label{eq:Sigmaz}
D_t \Sigma_{zz} 
&
+ 
2(1-\zeta)\Sigma_{zz} (\partial_xU + \partial_y V)
= (-2\nu(\partial_xU+\partial_yV)-\Sigma_{zz}) / \lambda 
\end{aligned}
\right.
$$
which is 
Johnson-Segalman's 3D model\footnote{
  Recall that Johnson-Segalman's model uses the family of Gordon-Schowalter \emph{objective} derivatives 
  ($\zeta=0$ is the upper-convected derivative, $\zeta=1$ is Jauman derivative, $\zeta=2$ is the lower-convected derivative)
  to cover a number of standard multidimensional generalizations of the Maxwell-fluid model:
  $\zeta=0$ in Oldroyd-B model, $\zeta=1$ in the co-rotational model and $\zeta=2$ in Oldroyd-A model
  (with an additional purely-viscous stress characterized by a retardation-time). 
} \cite{JOHNSON1977255} with slip parameter $\zeta\in[0,2]$ 
for some particular \emph{incompressible} velocity fields (those without vertical shear). 
%

The quasilinear system (\ref{eq:SVH}--\ref{eq:SVHUV}--\ref{eq:Sigmah}--\ref{eq:Sigmazz}) for $H$, $\bU$ 
and 
$\bC:=\lambda\bSigma/\nu+\bI$ is 
rotation-invariant. Then, it is easily computed that: 
\begin{proposition} 
Only the slip-parameter value $\zeta=0$ (i.e. the Upper-Convected 
case) 
ensures hyperbolicity of 
(\ref{eq:SVH}--\ref{eq:SVHUV}--\ref{eq:Sigmah}--\ref{eq:Sigmazz})
under the 
strain-free constraints: $H>0$, $C_{zz}>0$
and $\bC_h$ be a symmetric positive definite tensor (thus also $\bC$).
\end{proposition}

The proof follows after computing the eigenvalues of the jacobian in a 1D projection of the system (\ref{eq:SVH}--\ref{eq:SVHUV}--\ref{eq:Sigmah}--\ref{eq:Sigmazz}) like
\beq
\label{eq:svucm1D}
\left\lbrace
\begin{aligned}
\partial_{t} \ro  + \partial_{x}( \ro \ux ) & = 0
\\
\partial_{t} \ux  + \ux\partial_{x} \ux &
+ g \partial_x \ro 
- G \left(  
  (\sx -\sz)/\ro \partial_x \ro
+ \partial_x (\sx -\sz)
  \right)
= 0
\\
\partial_{t}      \uy  
+ \ux\partial_{x} \uy 
&
- G \left( 
  \sxy/\ro \partial_x \ro
+ \partial_x\sxy
 \right)
= 0
\\
\partial_{t}      \sx  + \ux\partial_{x} \sx 
&
- \left( 
2 (1-\zeta)\sx \partial_x \ux 
-\zeta\partial_x\uy 
\right) 
= 0
\\
\partial_{t}      \sy  
+ \ux\partial_{x} \sy 
&
- 
\sxy 
(2-\zeta)\partial_x\uy 
= 0
\\
\partial_{t}      \sxy  
+ \ux\partial_{x} \sxy 
&
- 
\left( (1-\zeta/2)\sx-\zeta/2)\sy \right) \partial_x\uy 
- (1-\zeta) \sxy \partial_x\ux 
= 0
\\
\partial_{t}      \sz  + \ux\partial_{x} \sz 
&
+ 2(1-\zeta)\sz \partial_x\ux 
= 0
\end{aligned}
\right.
\eeq
similarly to the proof in \cite{EDWARDS1990411} for a similar system written in stress variables $\bSigma$ when $\zeta=0$
(though without vertical stress and strain components, which allow here mass preservation).
Denoting $\Delta = 2gh + G\left( 2(3-2\zeta)\sz+\zeta\sy-3\zeta\sx \right)$, four eigenvalues read
$$
\ux\pm\frac12\sqrt{ 
 \Delta 
 + G\left( (4-2\zeta)\sx-2\zeta\sy \right) 
\pm\sqrt{ 
\Delta^2
+ G^2(4\zeta\sxy)^2
}
}
$$
and are real if, and only if, the following strain-parametrized inequality holds:
$$
G^2(4\zeta\sxy)^2 \le 2G\Delta\left( (4-2\zeta)\sx-2\zeta\sy \right) + G^2\left( (4-2\zeta)\sx-2\zeta\sy \right)^2
$$
where, however, the strain values $\sx,\sy,\sxy,\sz$ vanish 
when $\zeta=0$.
We therefore consider only 
(\ref{eq:SVH}--\ref{eq:SVHUV}--\ref{eq:Sigmah}--\ref{eq:Sigmazz}) when $\zeta=0$,
where hyperbolicity is ensured with eigenvalues $u\pm\sqrt{gh+3G\sz+G\sx}$, $u\pm\sqrt{G\sx}$ and $u$
(with multiplicity 3) under the physcially-natural constraints $h\ge0,\sz\ge0,\sx\ge0$.

\mycomment{
DOMAIN INVARIANT in hoff def ? \cite{boyaval-lelievre-mangoubi-2009}
}

Indeed, at each point of the flow,
$\bC$ 
can be interpreted as the expectation of $\bR\bR^T$ i.e. a \emph{covariance} matrix for a stochastic 
material vector field 
$\bR$ 
that elastically deforms following the overdamped Langevin equation
$$
d \bR = \left( - (\bu\cdot\grad) \bR + (\bL-\zeta\bD) \bR - \frac1\lambda \bR \right)dt + \frac1{2\sqrt\lambda} d\bB_t
$$ 
under a Brownian 
field 
$(\bB_t)$ \cite{ottinger-1996}.
Alternatively, the stress $\Strs=G(\bC-\bI)$ is also reminiscent of a linear elastic (Hookean) material 
with elastic 
modulus $G=\nu/\lambda$,
where 
the time-rate of $\bC$ would be that of a \emph{left} Cauchy-Green 
deformation tensor $\bF\bF^{T}$
associated with a deformation gradient $\bF$ of time-rate $D_t\bF-\gbu\bF$. 
This is easily seen on rewriting Saint-Venant-Upper-Convected-Maxwell (SVUCM) model with constitutive equations:
\begin{align}
\label{eq:Ch}
D_t \bC_h - \bL_h\bC_h - \bC_h \bL_h^T 
& = (\bI-\bC_h)/\lambda \,,
\\
\label{eq:Cz}
D_t C_{zz} 
 + 2 
 C_{zz} \div\bU & = (1-C_{zz}) / \lambda \,.
\end{align}


As expected, the SVUCM model defined by the hyperbolic quasilinear system (\ref{eq:SVH}--\ref{eq:SVHUV}--\ref{eq:Ch}--\ref{eq:Cz})
on the domain $H,C_{zz},\bC_h=\bC_h^T>0$ 
preserves mass.
Moreover, it 
formally satisfies the additional conservation law \eqref{eq:energyconservation} 
with dissipation 
$D:=G(\tr\bC+\tr\bC^{-1}-2\tr\bI)/(2\lambda)>0$
for the Helmholtz free-energy
\begin{equation}
\label{eq:energy} 
E = 
\left( |\bU|^2 + gH + G \tr(\bC-\ln\bC-\bI) \right)/2 \,,
\end{equation}
which suggests the following second thermodynamics principle formulation 
\begin{multline}
\label{eq:SVHE}
\partial_t(HE) + \partial_x\left( HEU + H(P+\Sigma_{zz}-\Sigma_{xx})U - H\Sigma_{xy}V \right) 
\\
+ \partial_y\left( HEV - H\Sigma_{yx}U + H(P+\Sigma_{zz}-\Sigma_{yy})V \right) \le -kH |\bU|^2 - HD 
\end{multline}
on the admissibility domain  $H,C_{zz},\tr\bC_h,\det\bC_h>0$ 
(where we use $\tr\ln\bC_h=\ln\det\bC_h$).
Note that the free energy reads $E=(U^2+V^2)/2+E_H+E_\Sigma$ 
with Saint-Venant's 
potential energy $E_H=gH/2$
plus a Hookean contribution 
$$ 
E_\Sigma(\bC) \equiv E_{\Sigma_h}+E_{\Sigma_{zz}}= 
{G} \left( \tr(\Ecal(\bC_h))+\Ecal(C_{zz}) \right)/2
$$ 
where the function $\Ecal(x):=x-\ln x-1$ is convex in $x>0$, but 
the system (\ref{eq:SVH}--\ref{eq:SVHUV}--\ref{eq:Sigmah}--\ref{eq:Sigmazz})
is not 
in purely conservative form
so 
one cannot straightforwardly apply Godunov-Mock theorem to show (symmetric) hyperbolicity of the system.
However, $E$ is useful to hyperbolicity: formally, the physically-admissible solutions to Cauchy problems
that satisfy \eqref{eq:SVHE} (and have admissible initial value in the hyperbolicity domain) 
preserve $C_{zz},\tr\bC_h,\det\bC_h>0$.

\medskip

However, note that there is a discrepancy in between the strain measure 
$(\bI-\bF^{-T}\bF^{-1})/2$ 
naturally associated with the \emph{left} Cauchy-Green (specific) deformation tensor $\bF\bF^{T}$ 
(usually termed Euler-Almansi) 
and the 
definition of the 
elastic stress $\Strs$.
That is why, we also propose another model where (\ref{eq:SVH}--\ref{eq:SVHUV}) is coupled through
$\bSigma_h=G(\bI-\bC_h)$, $\Sigma_{zz}=G(1-C_{zz})$ 
to
\begin{align}
\label{eq:Ch3}
D_t \bC_h + \bL_h\bC_h + \bC_h \bL_h^T 
& = (\bI-\bC_h)/\lambda \,,
\\
\label{eq:Cz3}
D_t C_{zz} - 2
C_{zz} \div\bU & = (1-C_{zz}) / \lambda \,.
\end{align}
Indeed, 
using for the time-rate of $\bC$ 
that of 
the \emph{inverse} left Cauchy-Green deformation tensor (sometimes termed Finger deformation tensor)
seems more 
consistent with the definition of (linear) elastic stresses $\bSigma$ as a function of 
a (Euler-Almansi) strain measure. 
Note that the time-rate 
tensor in (\ref{eq:Ch3}--\ref{eq:Cz3}) 
is Galilean-invariant
and induces an objective Maxwell-fluid law. 
It has been used previously in another context where 2D \emph{weakly-sheared} flows also appear
\cite{teshukov-2007,Teshukov2007,hal-01529497} 
and this is why we call this a Saint-Venant-Teshukov-Maxwell (SVTM) model.
%
The SVTM model is 
rotation-invariant. 
Then, a 1D projection 
like 
\beq
\label{eq:system1D}
\left\lbrace
\begin{aligned}
\partial_{t}    \ro 
+ \partial_{x}( \ro \ux ) 
& = 0
\\
\partial_{t}      \ux  
+ \ux\partial_{x} \ux 
+ g \partial_x \ro 
&
+ G \left(  
  (\sx -\sz)/\ro \partial_x \ro
+ \partial_x (\sx -\sz)
  \right)
= 0
\\
\partial_{t}      \uy  
+ \ux\partial_{x} \uy 
&
+ G \left( 
  \sxy/\ro \partial_x \ro
+ \partial_x\sxy
 \right)
= 0
\\
\partial_{t}      \sx  
+ \ux\partial_{x} \sx 
&
+ 2 \sx\partial_x \ux 
= 0
\\
\partial_{t}      \sy  
+ \ux\partial_{x} \sy 
&
+ 2 \sxy\partial_x\uy 
= 0
\\
\partial_{t}      \sxy  
+ \ux\partial_{x} \sxy 
&
+ \sxy \partial_x\ux 
+ 
\sx\partial_x\uy
= 0
\\
\partial_{t}      \sz  
+ \ux\partial_{x} \sz 
&
- 2 \sz 
\partial_x\ux
= 0
\end{aligned}
\right.
\eeq
again suffices (using lower-case notations for 1D systems)  
to see that the 
system is hyperbolic (strictly) 
on the same admissibility domain as 
SVUCM:
$$
\{ \ro,\sx,\sy,\sz>0\,;\quad (\ux,\uy)\in\R^2\,;\quad \sx\sy > \sxy^2 \} \,.
$$
In 
system 
\eqref{eq:system1D} 
for $Q=(h,u,v,\sx,\sxy,\sy,\sz)$, the jacobian 
has eigenvalues:
$$
\lambda_0 = \ux \text{ (with multiplicity 3),}
\
\lambda_{1\pm} = \ux \pm \sqrt{G\sx}
\,,\
\lambda_{2\pm} = \ux \pm \sqrt{g\ro+G(3\sx+\sz)} \,.
$$

The SVTM model satisfies the same mass and energy conservation laws as SVUCM,
but the laws (and the interpretations) of the dependent variable $\bC$ are not the same,
as we shall 
observe in numerical illustrations of Section~\ref{sec:num}.


\medskip

Finally, to sum up,
we can now propose two hyperbolic quasilinear systems with zero retardation time
as models for viscoelastic incompressible shallow free-surface 3D flows under gravity, 
which generalizes Saint-Venant's approach from a Newtonian fluid 
to Maxwell-type fluids. 

The first 2D model is one particular case of the depth-averaged 
model formally derived in \cite{bouchut-boyaval-2015}.
It bears similarities with the hyperbolic model of \cite{PHELAN1989197,EDWARDS1990411,Guaily2010158,Guaily2011258}
for compressible viscoelastic flows.
In comparison, note that our model applies to incompressible flows and \emph{conserves mass},
although it remains mostly 2D.
The improvement is obtained thanks to the introduction of a free-surface and assuming a stratified flow
(i.e. vertical acceleration is neglected and a vertical velocity can be reconstructed from mass conservation),
with a total (specific) pressure $P+\Sigma_{zz}$ decomposing into a hydrostatic component $P$ plus an elastic component $\Sigma_{zz}$ hence it holds
$$
D_t (P+\Sigma_{zz}) + \left( P + \frac{2}G \Sigma_{zz} + 2 \right) \div\bU = 0
$$
to be compared with 
$
D_t (P+\Sigma_{zz}) + \gamma \left( P + \Sigma_{zz} \right) \div\bU = 0
$
in \cite{PHELAN1989197,EDWARDS1990411,Guaily2010158,Guaily2011258}.

The second 2D model is an alternative to the first one
which also conserves mass and is compatible with a formulation of thermodynamics principles.
It has a more clear stress-strain relationship than the first model.
It has been used previously to model eddies-microstructures in (weakly-sheared) turbulent flows. 

Note that the two models contain similar closed subsystems that fully describe the 1D motions,
for $q_1=(\ro,\ux,\sz,\sx)$ in solutions invariant e.g. by translation along $\be_y$. 
In SVUCM case, this is the 1D model that we proposed in \cite{bouchut-boyaval-2013},
and whose Riemann initial-value problem has been completely solved in \cite{2016arXiv161108491B}.
In SVTM case, the 1D system reads similarly except that $\sx$ and $\sz$ exchange their roles.
This is natural insofar as their interpretation as strain components has been inverted~!
One should 
keep in mind the inverted interpretations of the variable $\bC$ as deformations (in SVUCM) or as inverse deformations (in SVTM),
in particular to read numerical results in Section~\ref{sec:num}.


\section{Computational schemes for 
Cauchy problems} 
\label{sec:FV}

\subsection{A framework for 
Riemann-based Finite-Volume 
schemes} 
\label{sec:general}

With a view to 
numerically simulating time evolutions of flows governed by our models SVUCM and SVTM, 
we consider now computable 
approximations of solutions to Cauchy problems for the hyperbolic quasilinear systems.

We 
study carefully the SVTM case (\ref{eq:SVH}--\ref{eq:SVHUV}--\ref{eq:Ch3}--\ref{eq:Cz3})
with $P=gH/2$ 
and $\bSigma_h=G(\bI-\bC_h)$, $\Sigma_{zz}=G(1-C_{zz})$ 
\begin{align*}
& \partial_t H + \div( H \bU ) = 0
\\
& \partial_t(H\bU) + \div(H\bU\otimes\bU + gH^2/2 - GH C_{zz} + GH\bC_h ) = -kH\bU
\\
& \partial_t \bC_h + \bU\cdot\grad \bC_h + 
+ \bL_h\bC_h + \bC_h \bL_h^T 
= (\bI-\bC_h)/\lambda
\\
& \partial_t C_{zz} + \bU\cdot\grad C_{zz} + 
- 2 
C_{zz} \div\bU = (1-C_{zz}) / \lambda
\end{align*}
first, then a similar numerical scheme will follow for SVUCM (\ref{eq:SVH}--\ref{eq:SVHUV}--\ref{eq:Ch}--\ref{eq:Cz})
\begin{align*}
& \partial_t H + \div( H \bU ) = 0
\\
& \partial_t(H\bU) + \div(H\bU\otimes\bU + gH^2/2 + GH C_{zz} - GH\bC_h ) = -kH\bU
\\
& \partial_t \bC_h + \bU\cdot\grad \bC_h + 
- \bL_h\bC_h - \bC_h \bL_h^T 
= (\bI-\bC_h)/\lambda
\\ 
& \partial_t C_{zz} + \bU\cdot\grad C_{zz} + 
2 
C_{zz} \div\bU = (1-C_{zz}) / \lambda
\end{align*}
with $\bSigma_h=G(\bC_h-\bI)$, $\Sigma_{zz}=G(C_{zz}-1)$ 
as explained in Section~\ref{sec:scheme},
on noting that both models have the same admissibility domain for the variable $\bC$, the same equilibrium, 
and that smooth solutions 
satisfy the same additional conservation law \eqref{eq:energyconservation} in both cases
(though we recall that the physical interpretations of the variable $\bC$ are inversed one another in between the two models).



We are not aware of a general theory that ensures the existence 
of a \emph{unique} solution to Cauchy initial-value problems
for \emph{multidimensional nonlinear} 
hyperbolic systems like SVTM (or SVUCM) model.
%
For instance, the 2D 
Saint-Venant 
system of conservation laws (formally reached 
in the limit $G\to0$),
which is symmetric and strictly hyperbolic, 
may have multiple \emph{entropy solutions}
\cite{chiodaroli-delellis-kreml-2015}.
%


\mycomment{ 
Pragmatically i.e. for numerical purposes, 
we shall heuristically neglect the tangential propagation of information (waves) in comparison with transverse at each cell interface,
as usual for \emph{rotation-invariant} hyperbolic systems (where information propagates at finite speed). Then
}

Yet, full-2D admissible solutions 
can be approximated \emph{stably} with a \emph{Riemann-based} Finite-Volume (FV) method,
using a piecewise-constant discretization 
of the unknown fields $H,\bU,\bC$ on the cells $V_i$ ($i\in \NN$) of a 
mesh \cite{leveque-2002}. 

In Riemann-based FV methods, 
one only needs numerical fluxes built from 1D Riemann 
problems
(formula \eqref{eq:conservation} below). 
Riemann problems have simple (self-similar) solutions in 1D,
and they can typically be shown \emph{well-posed} (with \emph{entropy solutions})
for strictly hyperbolic systems of genuinely-nonlinear conservation laws,
either by the vanishing-viscosity method or (equivalently) by the front-tracking method \cite{bianchini-bressan-2005}.

Recall however that 
our models are 
not full sets of conservation laws here.
Defining univoque solutions to Riemann problems for the 1D SVTM system \eqref{eq:system1D} 
may then be an issue because of non-conservative products,
see e.g. \cite{berthon-coquel-lefloch-2011} for a discussion of the problem and solutions.
One solution 
has been used in \cite{bouchut-boyaval-2013,2016arXiv161108491B} 
to univoquely define (approximate) 1D Riemann solutions for the closed $4\times 4$ subsystem in $(h,u\sx,\sz)$,
because the non-conservative variables correspond to linearly degenerate fields and
the mathematical entropy $HE$ remains convex on the Hugoniot loci.
But that solution method does not straightforwardly apply anymore here 
insofar as a non-conservative variable 
is not linearly degenerate.

Pragmatically though (i.e. for numerical purposes),
we only need univoque \emph{numerical fluxes}.
That is the reason why we are naturally led to using only \emph{approximate} (simple) 1D Riemann solvers
which use only simple solutions (with only a few degrees of freedom)
that are explicitly constrained to satisfy only the few consistency conditions identified so far,
see e.g. \cite{bouchut-2004}.

We will construct in Section~\ref{sec:oned}
\emph{univoque} physically-admissible simple (approximate) solutions to 1D Riemann problems (without source term),
thanks to a \emph{frame-invariant} Riemann solver ensuring a discrete version of our formulation of thermodynamics principles
(for the homogeneous systems with $k=0$).

Now, a framework exists such that a ``good'' simple Riemann solver ensures the consistency of Riemann-based FV approximations
with solutions to the homogeneous system.
Moreover we also have to handle a source term.
So, before defining precisely solutions to 1D Riemann homogeneous problems,
let us first recall that FV framework from \cite{bouchut-2004} and how to deal with the source term.




Given a tesselation of $\R^2$ using polygonal cells $V_i$, $i\in\NN$,
we write $q_h=\sum_i q_i(t) 1_{V_i}$ the FV discretization of a $7$-dimensional state vector $q$ for SVTM model.
%
(See below why we do not choose $Q\equiv(\ro,\ux,\uy,\sx,\sy,\sxy,\sz)$ as discretization variable $q$ for \eqref{eq:system1D}.)
We would like to solve a Cauchy problem for $q_h(t)$, $t\ge0$ given some admissible $q_h^0$ at $t=0\equiv t^0$.

\subsubsection{Splitting the time integration of FV approximations}

We use a splitting method to divide time-integration at time $t^n=\sum_{k=0}^{n-1} \tau^k$, $n=1\ldots N$ ($\tau^k>0$) 
into two sub-steps: 
first, the 
differential terms are integrated forward in time; 
second, the source terms are integrated backward in time. 

For time step 
$[t^n,t^{n+1})$, we start with $q^n_h\approx q(t^n)$ and first have to compute the 
solution $q^{n+1,-}_h=\sum_i q^{n+1,-}_i 1_{V_i}$ at $t^{n+1,-}$
of the Cauchy problem for the \emph{homogeneous} 
SVTM model 
on $[t^n,t^{n+1})$ with a forward-Euler time 
scheme.

Denoting $\bn_{i\to j}$ the unit normal at a face $\Gamma_{ij}\equiv\overline{V_i}\cap\overline{V_j}$ oriented from $V_i$ to $V_j$,
we compute 
$q^{n+1,-}_h$ 
with a FV approximation
on each control volume 
$V_i$
\begin{equation}
\label{eq:conservation}
q_i^{n+1,-} = q_i^n - \tau^n \sum_{ \Gamma_{ij}\equiv\overline{V_i}\cap\overline{V_j} \neq \emptyset} 
\frac{|\Gamma_{ij}|}{|V_i|} \bF_{i\to j}(q_i^n,q_j^n;\bn_{i\to j}) 
\end{equation} 
that is consistent with 
a time-integration of our \emph{frame-invariant model without source term}.
In particular, discrete conservation laws hold for 
$H_i$ and $H_i\bU_i$ with consistent numerical fluxes $\bF_{i\to j}^{H,H\bU}(q_i,q_j;\bn_{i\to j})$
for the $H$ and $H\bU$ 
components 
of $q$: the latter
should equal the 
normal flux components $\bF^{H,H\bU}(q)\cdot\bn_{i\to j}$ 
when $q=q_i=q_j$,
and satisfy
$ \bF_{i\to j}^{H,H\bU}(q_i,q_j;\bn_{i\to j}) + \bF_{j\to i}^{H,H\bU}(q_j,q_i;\bn_{j\to i}) = 0 $. 
For the other components, they are a priori solutions to non-conservative equations and
it is still not clear which consistency conditions should be satisfied for them at this stage\dots 
except 
i) the additional conservation law for $HE$, 
which we require here as the (Clausius-Duhem) \emph{inequality} \eqref{eq:energydissipation}, and
ii) the \emph{frame-invariance},
which is preserved if we use a frame-invariant 1D Riemann solver. 
%
%
%

We will define 
\emph{only in the next Section \ref{sec:oned}}
such a numerical flux of the form
\begin{equation}
\label{eq:rotationinvariantflux} 
\bF_{i\to j}(q_i,q_j;\bn_{i\to j})=\Oij\tilde\bF(\Oij^{-1}q_i,\Oij^{-1}q_j)
\end{equation}
that preserves the frame invariance 
and uses a 
flux $\tilde\bF$ consistent for 1D Riemann problems, 
where $\Oij^{-1}q$ denotes a vector 
state variable with components in the local basis $(\bn_{i\to j},\bn_{i\to j}^\perp)$ 
computed from the FV vector 
state $q$ defined in a fixed reference frame $(\be_x,\be_y)$
and $\Oij$ denotes the inverse operation.

Still, at this point, we can already recall 
standard \emph{stability properties} that shall be transfered to the FV approximation
if the numerical flux is well chosen.

If the 1D Riemann solver $\tilde\bF$ satisfies some stability constraints,
then one can ensure stability properties for the FV approximation $q_h^{n+1,-}$ of the homogeneous hyperbolic system 
like the preservation of admissible domains, or discrete second-principle formulations (under conditions, like CFL: see below).


\subsubsection{Sub-step 1: requiring CFL and entropy-consistency 
conditions} 

Let us consider numerical fluxes in (\ref{eq:conservation}--\ref{eq:rotationinvariantflux}) of the form
\begin{multline}
\label{eq:numfluxHLL}
\tilde\bF_{i\to j}(\Oij^{-1}q_i,\Oij^{-1}q_j)
\\
= \Oij^{-1}\bF(q_i)\bn_{i\to j}-\int_{-\infty}^0 \left( R(\xi,\Oij^{-1}q_i,\Oij^{-1}q_j)-\Oij^{-1}q_i\right)d\xi
\end{multline}
using a simple 1D Riemann solver $R$ 
with finite maximal 
wavespeed $s(q_l,q_r)>0$, 
which is 
consistent with 
conservation laws 
(when relevant: 
recall here only components $H$ and $H\bU$ in $q$ would satisfy conservation laws, a priori).
If the 
Riemann solver preserves some physically-meaningful domain like (invariant) sets convex in the discretization variable $q$, 
so will (\ref{eq:conservation}--\ref{eq:rotationinvariantflux}--\ref{eq:numfluxHLL}) do on noting
\begin{multline}
\label{eq:cvxcombin1}
q_i^{n+1,-} = q_i^n \left( 1 - \tau^n \sum_j \frac{|\Gamma_{ij}|s(\Oij^{-1}q_i^n,\Oij^{-1}q_j^n) }{|V_i|} \right) 
\\
+ \sum_j \frac{|\Gamma_{ij}|}{|V_i|} \int_{-s(\Oij^{-1}q_i^n,\Oij^{-1}q_j^n) \tau^n }^0 \Oij R(x/\tau^n,\Oij^{-1}q_i^n,\Oij^{-1}q_j^n) dx
\end{multline}
is a convex combination, under the CFL condition
\beq
\label{eq:CFLsoft}
\forall i \quad \tau^n \sum_{j} \frac{|\Gamma_{ij}|s(\Oij^{-1}q_i^n,\Oij^{-1}q_j^n) }{|V_i|} \le 1 \,.
\eeq

\begin{proposition}
Under CFL condition \eqref{eq:CFLsoft}, the FV time-integration formula \eqref{eq:conservation}
with flux \eqref{eq:rotationinvariantflux} based on the 1D Riemann solver \eqref{eq:numfluxHLL}
preserves state sets that are convex in the discretization variable $q$. 
\end{proposition}

Moreover, a more stringent CFL condition with $s_i^n := \max_j s(\Oij^{-1}q_i^n,\Oij^{-1}q_j^n)$
\beq
\label{eq:CFL}
\dt^n s_i^n \le \left( \sum_j \frac{|\Gamma_{ij}|}{|V_i|} \right)^{-1}
\eeq
may allow one to use Jensen inequality with the convex combination 
\beq
\label{eq:cvxcombin2}
q_i^{n+1,-} 
= \sum_j \frac{|\Gamma_{ij}|}{|V_i|} 
 \int_{-\left(\sum_j \frac{|\Gamma_{ij}|}{|V_i|}\right)^{-1}}^0 \Oij R(x/\tau^n,\Oij^{-1}q_i^n,\Oij^{-1}q_j^n) dx \,,
\eeq
so as to ensure \emph{not only the decay of functionals $S(q)$ convex in $q$} 
(already ensured by \eqref{eq:cvxcombin1} when the Riemann solver preserves convex invariant sets)
but also a consistent discrete version of a \emph{second-principle formulation} like
\begin{equation}
\label{eq:secondprinc}
\partial_t S(q) + \div 
\bG(q) \le 0
\end{equation}
whenever \eqref{eq:secondprinc} holds for some ``mathematical entropy'' $S(q)$
that is a convex function of the \emph{Galilean invariants} of the state $q$ 
with 
``entropy-flux'' 
$\bG(q)$,
\emph{under additional entropy-consistency conditions} on the 
solver (see \eqref{solversufficientcondition} below). 
%
%
The consistency of numerical approximations with second-principle inequalities is essential
to possibly converge to physically admissible \emph{entropy solutions} 
\cite{dafermos-2000}. 

\medskip

Now, recall the SVTM model (\ref{eq:SVH}--\ref{eq:SVHUV}--\ref{eq:Ch3}--\ref{eq:Cz3}) is complemented by \eqref{eq:SVHE} i.e.
\begin{multline}
\label{eq:SVTMHE}
\partial_t(HE) + \partial_x\left( HEU + H(P+C_{xx}-C_{zz})U + HC_{xy}V \right) 
\\
+ \partial_y\left( HEV + HC_{yx}U + H(P+C_{yy}-C_{zz})V \right) \le -kH |\bU|^2 - HD 
\end{multline} 
with $D=-G(4-\tr\bC_h-\tr\bC_h^{-1}+2-C_{zz}-C_{zz}^{-1})/(2\lambda)\ge0$ for
$E=(U^2+V^2)/2+E_H+E_\Sigma$ 
where 
$E_H=gH/2$, 
$ 
E_\Sigma(\bC) = 
{G} \left( \tr(\Ecal(\bC_h))+\Ecal(C_{zz}) \right)/2 \,,
$
and $\Ecal(x):=x-\ln x-1$. 
%
%
%
Then, we consider $S=HE$ as a mathematical entropy for 
the SVTM system without source term $k=0=D$ (in the first sub-step of our time-splitting scheme).
The inequality \eqref{eq:SVTMHE} shows that a second-principle formulation \eqref{eq:secondprinc} holds with an entropy flux 
$\bG=HE\bU+H(P+\Sigma_{zz})\bU-H\bSigma_h\cdot\bU$.

\medskip

So first, with a view to using Jensen inequality with \eqref{eq:cvxcombin2},
we can already choose a FV discretization variable 
$q$ with a convex admissible set such that $S\equiv HE$ is convex in $q$ at this stage.
We propose 
\beq
\label{convexvar}
 q:=\left(H,HU,HV,HC_{xx},HC_{yy},HC_{xy}/\sqrt{C_{xx}C_{yy}},HC_{zz}\right)
\eeq
which has a convex admissible domain $q_1,q_4,q_5,q_7>0,|q_6|<q_1$
and such that $S/H\equiv E$ is convex  
in $(H^{-1},U,V,C_{xx},C_{yy},C_{xy}/\sqrt{C_{xx}C_{yy}},C_{zz})$,
recall \cite[Lemma 1.4]{bouchut-2004}.

Note however that the Riemann problems at interfaces $\Gamma_{ij}$ will be solved 
for another variable $\tilde q$, 
see Section \ref{sec:oned}.
So, 
let us stress again that
the operators $\Oij^{-1}$ are \emph{nonlinear}\footnote{
  The operators would degenerate as linear functions of the vector state 
  only if we used the same representation at each step of the algorithm.
  Then, the operations would simply consist in a -- linear -- change of frame
  with coefficients \emph{quadratic} in the components of 
  the rotation matrix
  $$
  \bOmega_{ij}
  = \begin{pmatrix}
    \bn_{i\to j}\cdot\be_x & \bn_{i\to j}^\perp\cdot\be_x
    \\
    \bn_{i\to j}\cdot\be_y & \bn_{i\to j}^\perp\cdot\be_y
    \end{pmatrix}
  $$
  here, because $\bC_h$ are components of a (2-contravariant) 2-tensor variable.
} 
(functions of the 
vector representation $q$ of the state in a fixed Cartesian reference frame $\be_x,\be_y$),
as well as the inverse operators $\Oij$ 
(functions of the vector representation $\tilde q$ of the state in a local Cartesian basis $\bn_{i\to j},\bn_{i\to j}^\perp$).

\medskip

Next, to complete the first sub-step of our time-splitting scheme (a homogeneous Riemann problem),
we propose to build a 1D Riemann solver that satisfies the following \emph{entropy-consistency condition}: 
there exists a \emph{conservative} discrete entropy-flux $ \tilde G(q_l,q_r)=-\tilde G(q_r,q_l)$, 
\emph{consistent} in FV 
sense with \eqref{eq:secondprinc} 
$$ \tilde G(\Oij^{-1}q,\Oij^{-1}q)=\bG(\Oij^{-1}q)\cdot\be_x=\bG(q)\cdot\bn_{i\to j} $$ 
for all admissible states $q$, such that for any admissible states $q_l,q_r$ there holds
\begin{multline}
\label{solversufficientcondition}
%
\bG(q_r)\cdot\be_x + \int\displaylimits_0^{+\infty} \left( S\left(R(\xi,q_l,q_r)\right) -S(q_r) \right)d\xi
\\
\le
\tilde G(q_l,q_r)
\le 
\bG(q_l)\cdot\be_x - \int\displaylimits_{-\infty}^0 \left( S\left(R(\xi,q_l,q_r)\right) -S(q_l) \right)d\xi \,.
\end{multline}
Indeed, using
$ 
\int\displaylimits_0^{+\infty} S\left(R(\xi,q_l,q_r)\right) = - \int\displaylimits_{-\infty}^0  S\left(R(\xi,q_r,q_l)\right) 
$
and \eqref{solversufficientcondition} we obtain
\begin{multline}
\label{eq:numfluxHLLintegralentropy1}
\int\displaylimits_{-s_i^n \tau^n}^0 S\left(\Oij R(x/\tau^n,\Oij^{-1}q_i^n,\Oij^{-1}q_j^n)\right)dx
=
\int\displaylimits_{-s_i^n \tau^n}^0 S\left(R(x/\tau^n,\Oij^{-1}q_i^n,\Oij^{-1}q_j^n)\right)dx
\\
\le \left(\tau^n s_i^n\right) \: S(q_i^n)
- \tau^n \left( \tilde G(\Oij^{-1}q_i^n,\Oij^{-1}q_j^n) - \bG(q_i^n)\cdot\bn_{i\to j} \right)
\end{multline}
and finally, with \eqref{eq:cvxcombin2}, a consistent discrete 
second-principle formulation:
\begin{multline}
\label{eq:discreteentropy}
S(q_i^{n+1,-}) - S(q_i^n) + \tau^n \sum_j \frac{|\Gamma_{ij}|}{|V_i|} \tilde G(\Oij^{-1}q_i^n,\Oij^{-1}q_j^n)
\\
\le 
\tau^n \sum_j \frac{|\Gamma_{ij}|}{|V_i|} \bG(q_i^n)\cdot\bn_{i\to j} \equiv 0 \,.
\end{multline}

\begin{proposition}
Under CFL condition \eqref{eq:CFL}, the FV time-integration formula \eqref{eq:conservation}
with flux \eqref{eq:rotationinvariantflux} based on a 1D Riemann solver \eqref{eq:numfluxHLL}
\begin{itemize}
 \item preserves convex state sets in the FV discretization variable $q$ 
and
 \item satisfies the discrete version \eqref{eq:discreteentropy} of the second-principle formulation \eqref{eq:SVTMHE} 
(with $k=0=D$) 
provided
\end{itemize}
\begin{enumerate}
 \item $S=HE$ is a convex function of $q$, and
 \item the entropy-consistency condition \eqref{solversufficientcondition} is satisfied for some $\tilde G$.
\end{enumerate}
\end{proposition}
We propose in next Section~\ref{sec:oned} a 
solver $R$ satisfying \eqref{solversufficientcondition} for some $\tilde G$.

\subsubsection{Sub-step 2: integrating source terms with Backward-Euler}

To complete the time-integration of our SVTM model, 
we standardly propose a second sub-step using the backward-Euler time-scheme for the source terms
\begin{align}
\label{cnvxcombinationu}
\bU(q_i^{n+1}) - \bU(q_i^{n+1,-}) & = -\tau^n k \bU(q_i^{n+1})
\\
\label{cnvxcombinationh}
\bC_h(q_i^{n+1}) - \bC_h(q_i^{n+1,-}) & = \frac{\tau^n}{\lambda} \left(\bI-\bC_h(q_i^{n+1})\right)
\\
\label{cnvxcombinationz}
C_{zz}(q_i^{n+1}) - C_{zz}(q_i^{n+1,-}) & = \frac{\tau^n}{\lambda} \left(1-C_{zz}(q_i^{n+1})\right)
\end{align}
with $H_i^{n+1} = H_i^{n+1,-}$ (mass is conserved),
the time step $\tau^n$ being given by the CFL condition \eqref{eq:CFL} of sub-step 1. 
An \emph{admissible} state $\bU(q_i^{n+1}),\bC_h(q_i^{n+1}),C_{zz}(q_i^{n+1})$ can be computed \emph{explicitly} here
as a convex combination of admissible states,
insofar as the source terms RHS in \eqref{cnvxcombinationu}, \eqref{cnvxcombinationh}, \eqref{cnvxcombinationz}
are all \emph{linear} in $\bU,\bC_h,C_{zz}$, of relaxation type.
(We only need \emph{explicit} nonlinear mappings for the variable change $(H,\bU,\bC_h,C_{zz})\leftrightarrow q$ \emph{cellwise}.)

A discrete version of \eqref{eq:SVTMHE} is satisfied (with $D\neq0$ like in \eqref{eq:SVTMHE})
\begin{multline}
\label{eq:discreteentropy2}
S(q_i^{n+1}) - S(q_i^n) + \tau^n \sum_j \frac{|\Gamma_{ij}|}{|V_i|} \tilde G(\Oij^{-1}q_i^n,\Oij^{-1}q_j^n)
\\
\le 
-\tau^n k |\bU_i^{n+1}|^2 -\tau^n D(q_i^{n+1})
\end{multline}
which can be shown from \eqref{eq:discreteentropy} and the convexity of 
$S=HE$ 
hence, for instance,
\begin{multline}
\label{eq:cvxtensor}
E_\Sigma(\bC_h(q_i^{n+1})) - E_\Sigma(\bC_h(q_i^{n+1,-})) \\
\le (\bI-\bC_h(q_i^{n+1})^{-1}):\left(\bC_h(q_i^{n+1}) - \bC_h(q_i^{n+1,-})\right)
\end{multline}
where $(\bI-\bC_h(q_i^{n+1})^{-1}):\left(\bI-\bC_h(q_i^{n+1})\right)\le 0$ is one term of the dissipation $D(q_i^{n+1})$
(for a detailed proof of \eqref{eq:cvxtensor}, see e.g. \cite[(2.7e) in Lemma 2.1]{barrett-boyaval-2011}).

\subsection{A 5-wave 
relaxation solver}
\label{sec:oned}

To follow the general framework presented in Section~\ref{sec:general} for Riemann-based FV discretizations,
let us start here the construction of an entropy-consistent \emph{simple} Riemann solver $R(\zeta,\Oij^{-1}q_i,\Oij^{-1}q_j)$.
Precisely, we explicitly define approximate 
solutions to 1D Riemann problems for the quasilinear (non-conservative) SVTM 1D system \eqref{eq:system1D}
that are piecewise-constant functions of the self-similarity variable $\zeta$ with finitely-many values.
We next show 
that the entropy-consistency condition \eqref{solversufficientcondition} can be satisfied,
so the Riemann solver confomrs with the general framework presented in Section~\ref{sec:general}.
We recall that the SVUCM system will be treated afterwards in Section \ref{sec:scheme}, with an approach similar to that for SVTM.

To start with, let us consider the system \eqref{eq:system1D} for 
the variable 
$$
\tilde q:=\left(h,hu,h^{-2}c_{xx},h^2c_{zz},h(c_{yy}-c_{xy}^2/c_{xx}),hc_{xy}/\sqrt{c_{xx}},hv\right)
$$
with left/right initial values $\Oij^{-1}q_i,\Oij^{-1}q_j$ computed in the local frame by 
$$
\begin{aligned}
c_{xx}&=\bC_h(q)\bn_{i\to j}\cdot\bn_{i\to j}
\\
c_{xy}&=\bC_h(q)\bn_{i\to j}^\perp\cdot\bn_{i\to j}=\bC_h\bn_{i\to j}^\perp\cdot\bn_{i\to j}
\\ 
c_{xy}&=\bC_h(q)\bn_{i\to j}^\perp\cdot\bn_{i\to j}^\perp
\end{aligned}
$$
as \emph{nonlinear} functions of the left/right values $q_i,q_j$. 

For consistency, we require that Riemann solutions preserve
$ \ro,\sx,\sz,\sy-\sxy^2/\sx>0 $ 
i.e. admissibility, 
and mimick the second-principle formulation 
\beq
\label{fulineq}
\partial_t\left( h(E_\parallel + E_\perp) \right) + \partial_x\left( h\ux(E_\parallel + hE_\perp + \ux P_\parallel + \uy P_\perp ) \right) \le 0
\eeq
for the free-energy $ h(E_\parallel + E_\perp) $ as mathematical entropy, with two terms
$$
E_\parallel  = \frac{u^2}2 + e_\parallel \;;\quad e_\parallel = \frac{gh}2 + \frac{G(\sx+\sz)}2 - \frac{G\ln(\sx\sz)}2 
$$
$$
E_\perp =  \frac{v^2}2 + e_\perp  \;;\quad e_\perp = \frac{G\sxy^2/\sx}2 + \frac{G(\sy-\sxy^2/\sx)}2 - \frac{G\ln(\sy-\sxy^2/\sx)}2
$$
that satisfy 
independent conservation laws 
with two pressure terms
$$
\partial_t\left( h E_\parallel \right) + \partial_x\left( h\ux_\parallel + \ux P_\parallel \right) = 0
\quad
\partial_t\left( h E_\perp \right) + \partial_x\left( \uy( h\uy E_\perp + \uy P_\perp \right) = 0
$$
$$
P_\parallel = -\partial_{h^{-1}}|_{h\sx^{-1/2},h\sz^{+1/2}} e_\parallel 
= \frac{gh^2}2 + G h(\sx-\sz)
\qquad
P_\perp 
= G\ro\sxy 
$$
in smooth evolutions (without discontinuities) so as to define a consistent numerical entropy-flux $\tilde G$
(our condition \eqref{solversufficientcondition} for admissible discretizations).
%
%
One issue is \emph{first} give a meaning to the non-conservative nonlinear terms.
Here, 
we straightforwardly devise a \emph{simple} Riemann solver 
that univoquely approximates 
admissible 1D solutions. 
%
%
%
Rewriting 
the 1D SVTM model in $\tilde q$:
\beq
\label{eq:system1D2}
\left\lbrace
\begin{aligned}
\partial_{t}    \ro 
+ \partial_{x}( \ro \ux ) 
& 
= 0
\\
\partial_{t}(  \ro\ux  )
+ \partial_{x}( \ro \ux^2 ) 
+ \partial_x( g \ro^2/2 ) 
+ \partial_x\left( G \ro(\sx -\sz) \right)
&
= 0
\\
\partial_{t}( \ro^{-2}\sx )  
+ \ux \partial_{x}( \ro^{-2}\sx ) 
&
= 0
\\
\partial_{t}(  \ro^2\sz  )
+ \ux\partial_{x}( \ro^2\sz  ) 
&
= 0
\\
\partial_{t}( \sy-\sxy^2/\sx  )
+ \ux\partial_{x}( \sy-\sxy^2/\sx  ) 
&
= 0
\\
\partial_{t}( \ro \uy )
+ \partial_{x}( \ro \ux \uy ) 
+ \partial_x\left( G \ro\sxy \right)
&
= 0
\\
\partial_{t}( \ro^{-1}\sxy )
+ \ux\partial_{x}( \ro^{-1}\sxy )
+ (\ro^{-2}\sx) h\partial_x\uy
&
= 0
\end{aligned}
\right.
\eeq
with a view to constructing 
a 
Riemann solver, 
our variable choice $\tilde q$ for the 1D Riemann problems obviously justifies: 
there is actually \emph{only one} non-conservative product (for the evolution of $c_{xy}h^{-1}$).


We next consider a \emph{5-wave 
simple solver} for the system \eqref{eq:system1D2} 
in Euler coordinates (i.e. which uses a Eulerian flow description),
after the standard 
transformation of a simple solver 
for \eqref{eq:system1D2} rewritten in Lagrange coordinates
\cite{gallice-2000} i.e.
\beq
\label{eq:system1D2lagrange}
\left\lbrace
\begin{aligned}
\partial_{t}   \ro^{-1} 
- \partial_{x} \ux 
& 
= 0
\\
\partial_{t} \ux  
+ \partial_x\left( gh^2/2 + Gh(\sx-\sz) \right) 
&
= 0
\\
\partial_{t}( \ro^{-2}\sx )  
&
= 0
\\
\partial_{t}(  \ro^2\sz  )
&
= 0
\\
\partial_{t}( \sy-\sxy^2/\sx  )
&
= 0
\\
\partial_{t} \uy
+ \partial_x\left( G \ro\sxy \right)
&
= 0
\\
\partial_{t}( \sxy/h )
+ \sx \partial_x\uy
&
= 0
\end{aligned}
\right.
\eeq
which suggests a 5-wave (Lagrangian) solver inspired by Suliciu's relaxation strategy of pressures
(a ``general'' strategy for fluids described in \cite{bouchut-2004}).
The additional conservations laws for \eqref{eq:system1D2lagrange} write
$$
\partial_t E_\parallel + \partial_x \left( \ux P_\parallel \right) 
= 0 =
\partial_t E_\perp + \partial_x \left( \uy P_\perp \right)
$$
for smooth solutions. 
Then, on noting smooth $ P_\parallel=\frac{gh^2}2 + Gh(\sx-\sz)$ satisfy
$$
\partial_t P_\parallel + c_\parallel^2 \partial_x \ux = 0
$$
with $c_\parallel^2=h^2(gh+G(3\sx+\sz))$, 
smooth $P_\perp=G\ro\sxy$ satisfy
$$
\partial_t P_\perp + b \partial_x \ux + c_\perp^2 \partial_x \uy = 0
$$
with $b=2G\ro^2\sxy$ and $c_\perp^2=G\ro^2\sx$, and smooth $\sxy/h$ 
satisfy
$$
\partial_{t} (\sxy/h) + a^2 \partial_x\uy = 0
$$
with $a^2=\sx$, we propose 
the following 
relaxed approximation of \eqref{eq:system1D2lagrange}
\beq
\label{eq:system1D2lagrangerelaxed}
\left\lbrace
\begin{aligned}
\partial_{t} \ro^{-1} - \partial_{x} \ux 
& 
= 0
\\
\partial_{t} \ux + \partial_x \pi_\parallel
&
= 0
\\
\partial_t \pi_\parallel + c_\parallel^2 \partial_x \ux
&
= 0
\\
\partial_{t}( \sxy/\ro ) + a^2 \partial_x\uy
&
= 0
\\
\partial_{t} \uy + \partial_x \pi_\perp 
&
= 0
\\
\partial_t \pi_\perp + b \partial_x \ux + c_\perp^2 \partial_x \uy
&
= 0
\\
\partial_{t}( \ro^{-2}\sx )  
&
= 0
\\
\partial_{t}(  \ro^2\sz  )
&
= 0
\\
\partial_{t}( \sy-\sxy^2/\sx  )
&
= 0
\end{aligned}
\right.
\eeq
which is a hyperbolic system with all fields linearly degenerate
and thus allows one to compute easily approximate Riemann solutions in Lagrange coordinates.

The general Riemann solution to \eqref{eq:system1D2lagrangerelaxed}
has 5 waves with speeds $-c_\parallel<-c_\perp<0<c_\perp<c_\parallel$
that can be ordered consistently 
with the definition 
of the relaxation parameters 
in smooth evolution cases.
Note the following relations
\beq
\label{eq:system1D2lagrangerelaxeddiag1}
\left\lbrace
\begin{aligned}
\partial_{t}(\pi_\parallel+c_\parallel \ux) + c_\parallel 
\partial_{x}(\pi_\parallel+c_\parallel \ux) 
& = 0
\\
\partial_{t}(\pi_\parallel-c_\parallel \ux) - c_\parallel 
\partial_{x}(\pi_\parallel-c_\parallel \ux) 
& = 0
\\
\partial_{t}(\pi_\parallel/c_\parallel^2+\ro^{-1})  
& = 0
\\
\partial_{t}( \ro^{-2}\sx )  
&
= 0
\\
\partial_{t}(  \ro^2\sz  )
&
= 0
\end{aligned}
\right.
\eeq
hold for the $q_1$-subsystem, already treated in \cite{bouchut-boyaval-2013} with same relaxation approach.

So, the Riemann solution 
has the following structure 
(in any variable $q$)
\begin{equation}
\label{eq:relaxationsolution}
\begin{cases}
 q_l 
 & x < -c_\parallel t \\
 q_l^\star 
 & -c_\parallel t < x < -c_\perp t \\
 q_l^\sharp 
 & -c_\perp t < x < 0 \\
 q_r^\sharp
 & 0 < x < c_\perp t \\
 q_r^\star 
 & c_\perp t < x < c_\parallel t \\
 q_r 
 & x > c_\parallel t
\end{cases}
\end{equation}
and $\tilde q_1$ can be 
explicited through the following analytical expressions
\begin{equation}
\label{solutiondetailed2}
\begin{array}{c}
u_l^\star=u_r^\star=u_l^\sharp=u_r^\sharp= 
u^*\equiv\frac{c_\parallel u_l+\pi_{\parallel,l}+c_\parallel u_r-\pi_{\parallel,r}}{2c_\parallel}
\\
\pi_{\parallel,l}^\star=\pi_{\parallel,r}^\star=\pi_{\parallel,l}^\sharp=\pi_{\parallel,r}^\sharp=
\pi^*\equiv\frac{\pi_{\parallel,l}+c_{\parallel}u_l+\pi_{\parallel,r}-c_{\parallel}u_r}{2}
\\
\frac{1}{h_l^*}=\frac{1}{h_l^\sharp}=\frac{1}{h_l}\left(1+\frac{c_{\parallel}(u_r-u_l)+\pi_{\parallel,l}-\pi_{\parallel,r}}{2c_{\parallel}^2/h_l}\right)
\\
\frac{1}{h_r^*}=\frac{1}{h_r^\sharp}=\frac{1}{h_r}\left(1+\frac{c_{\parallel}(u_r-u_l)+\pi_{\parallel,r}-\pi_{\parallel,l}}{2c_{\parallel}^2/h_r}\right)
\\
(\ro^{-2}\sx)_l^\star = (\ro^{-2}\sx)_l^\sharp = \ro^{-2}\sx
\quad
(\ro^{-2}\sx)_r^\star = (\ro^{-2}\sx)_r^\sharp = \ro^{-2}\sx
\\
(\ro^{2}\sz)_l^\star = (\ro^{2}\sz)_l^\sharp = \ro^{2}\sz
\quad
(\ro^{2}\sz)_r^\star = (\ro^{2}\sz)_r^\sharp = \ro^{2}\sz
\end{array}
\end{equation}
(i.e. $u,\pi_\parallel$ are weak Riemann invariants for $c_{-1},c_0,c_1$ waves,
$\pi_\parallel/c_\parallel^2+\ro^{-1},\ro^{-2}\sx,\ro^{2}\sz$ 
are weak Riemann invariants for $c_{-2},c_{-1},c_1,c_2$ waves,
$\pi_\parallel-c_\parallel \ux$ for $c_2$ and $\pi_\parallel+c_\parallel \ux$ for $c_{-2}$)
on recalling 
\cite{bouchut-boyaval-2013}.
%
Note also
\beq
\label{eq:sysstar}
\begin{aligned}
\partial_{t}(\pi_\perp+c_\parallel \uy) + c_\parallel 
\partial_{x}(\pi_\perp+c_\parallel \uy) 
+ \partial_x\left( b\ux+(c_\perp^2-c_\parallel^2)\uy \right)
& = 0
\\
\partial_{t}(\pi_\perp-c_\parallel \uy) - c_\parallel 
\partial_{x}(\pi_\perp-c_\parallel \uy) 
+ \partial_x\left( b\ux+(c_\perp^2-c_\parallel^2)\uy \right)
& = 0
\end{aligned}
\eeq
where $b\ux+(c_\perp^2-c_\parallel^2)\uy $ is a weak Riemann invariant for $c_{-2},c_2$ waves like $a^2\pi_\perp/c_\parallel^2-\sxy/h$.
Moreover, $\pi_\perp-c_\parallel \uy$ is a weak $c_2$ invariant
and 
$\pi_\perp+c_\parallel \uy$ a weak $c_{-2}$ invariant.
The following 
analytical expressions then also hold 
\begin{equation}
\label{solutiondetailed3}
\begin{array}{c}
v_l^\star
=v_l+\frac{b}{c_\perp^2-c_\parallel^2}(u_l-u^\star)
=v_l+\frac{b}{c_\perp^2-c_\parallel^2}\frac{c_\parallel u_l-\pi_{\parallel,l}-c_\parallel u_r+\pi_{\parallel,r}}{2c_\parallel}
\\
v_r^\star
=v_r+\frac{b}{c_\perp^2-c_\parallel^2}(u_r-u^\star)
=v_r+\frac{b}{c_\perp^2-c_\parallel^2}\frac{c_\parallel u_r+\pi_{\parallel,r}-c_\parallel u_l-\pi_{\parallel,l}}{2c_\parallel}
\\
\pi_{\perp,l}^\star
=c_\parallel(v_l-v_l^\star) + \pi_{\perp,l}
=\pi_{\perp,l} + \frac{b c_\parallel }{c_\perp^2-c_\parallel^2}
\frac{c_\parallel u_r-\pi_{\parallel,r}-c_\parallel u_l+\pi_{\parallel,l}}{2c_\parallel}
\\
\pi_{\perp,r}^\star
=c_\parallel(v_r^\star-v_r) + \pi_{\perp,r}
=\pi_{\perp,r} + \frac{b c_\parallel }{c_\perp^2-c_\parallel^2}
\frac{c_\parallel u_r+\pi_{\parallel,r}-c_\parallel u_l-\pi_{\parallel,l}}{2c_\parallel}
\\
\left(\sxy/h\right)_l^\star 
= \left(\sxy/h\right)_l 
+ \frac{a^2}{c_\parallel^2}(\pi_{\perp,l}^\star-\pi_{\perp,l})
= \left(\sxy/h\right)_l 
+ \frac{a^2 b / c_\parallel }{c_\perp^2-c_\parallel^2}
\frac{c_\parallel u_r-\pi_{\parallel,r}-c_\parallel u_l+\pi_{\parallel,l}}{2c_\parallel}
\\
\left(\sxy/h\right)_r^\star 
= \left(\sxy/h\right)_r
+ \frac{a^2}{c_\parallel^2}(\pi_{\perp,r}^\star-\pi_{\perp,r})
= \left(\sxy/h\right)_r 
+ \frac{a^2 b / c_\parallel }{c_\perp^2-c_\parallel^2}
\frac{c_\parallel u_r+\pi_{\parallel,r}-c_\parallel u_l-\pi_{\parallel,l}}{2c_\parallel}
\end{array}
\end{equation}
and, on noting
\beq
\label{eq:syssharp}
\begin{aligned}
\partial_{t}(\pi_\perp+c_\perp \uy) + c_\perp 
\partial_{x}(\pi_\perp+c_\perp \uy) + b \partial_x u
& = 0
\\
\partial_{t}(\pi_\perp-c_\perp \uy) - c_\perp 
\partial_{x}(\pi_\perp-c_\perp \uy) + b \partial_x u
& = 0
\\
\partial_{t}(a^2\pi_\perp-c_\perp^2\sxy/h) + a^2 b \partial_x u
& = 0
\end{aligned}
\eeq
where 
$u$ is a weak 
invariant for $c_{-1},c_0,c_1$, we get
\begin{equation}
\label{solutiondetailed4}
\begin{array}{c}
v_l^\sharp=v_r^\sharp=
v^\sharp\equiv\frac{c_\perp v_l^\star+\pi_{\perp,l}^\star+c_\perp v_r^\star-\pi_{\perp,r}^\star}{2c_\perp}
\\
\pi_{\perp,l}^\sharp=\pi_{\perp,r}^\sharp=
\pi_{\perp}^\sharp=\frac{\pi_{\perp,l}^\star+c_{\perp}v_l^\star+\pi_{\perp,r}^\star-c_{\perp}v_r^\star}{2}
\\
\left(\sxy/h\right)_r^\sharp 
= \left(\sxy/h\right)_r^\star 
+ \frac{a^2}{c_\perp^2}(\pi_{\perp}^\sharp-\pi_{\perp,r}^\star)
= \left(\sxy/h\right)_r^\star 
+ \frac{a^2}{c_\perp^2}\frac{\pi_{\perp,l}^\star+c_{\perp}v_l^\star-\pi_{\perp,r}^\star-c_{\perp}v_r^\star}{2}
\\
\left(\sxy/h\right)_l^\sharp 
= \left(\sxy/h\right)_l^\star 
+ \frac{a^2}{c_\perp^2}(\pi_{\perp}^\sharp-\pi_{\perp,l}^\star) 
= \left(\sxy/h\right)_l^\star 
+ \frac{a^2}{c_\perp^2}\frac{\pi_{\perp,r}^\star-c_{\perp}v_r^\star-\pi_{\perp,l}^\star+c_{\perp}v_l^\star}{2}
\end{array}
\end{equation}
which completes the 
expression of the Riemann solution of \eqref{eq:system1D2lagrangerelaxed} with, for $o=l/r$
$$
\left(\sy-\sxy^2/\sx\right)_o = \left(\sy-\sxy^2/\sx\right)_o^\star =\left(\sy-\sxy^2/\sx\right)_o^\sharp
\,.
$$

\smallskip

Now, we can propose a ``pseudo-relaxed'' Riemann solver 
for SVTM system 
\beq
\label{eq:system1Deulerrelaxed}
\left\lbrace
\begin{aligned}
\partial_{t} \ro + \partial_{x} (\ro\ux)
& 
= 0
\\
\partial_{t} (\ro\ux) + \partial_x (\ro\ux^2+\pi_\parallel)
&
= 0
\\
\partial_t (\ro\pi_\parallel) + \partial_x (\ro\ux\pi_\parallel) + c_\parallel^2 \partial_x \ux
&
= 0
\\
\partial_{t} \sxy + \partial_x (\ux\sxy)+ a^2 \ro\uy)
&
= 0
\\
\partial_{t}(\ro\uy) + \partial_x (\ro\ux\uy+\pi_\perp)
&
= 0
\\
\partial_t (\ro\pi_\perp) + \partial_x (\ro\ux\pi_\perp+ b \ux) + c_\perp^2 \partial_x \uy
&
= 0
\\
\partial_{t}( \ro^{-1}\sx ) + \partial_{x}( \ro^{-1}\sx\ux )  
&
= 0
\\
\partial_{t}(  \ro^3\sz  ) + \partial_{t}(  \ro^3\sz\ux  )
&
= 0
\\
\partial_{t}\left(  \ro( \sy-\sxy^2/\sx  )  \right) + \partial_{x}\left(  \ro( \sy-\sxy^2/\sx  )\ux  \right)
&
= 0
\end{aligned}
\right.
\eeq
in \emph{Euler coordinates}, 
recalling the 
transformation of a simple Riemann solver from Lagrange to Euler coordinates \cite{gallice-2000}.
Note that it has the same intermediate states as the solver in Lagrange coordinates, but different wave-speeds, namely:
\begin{align*}
\xi_{-2} & = u_l-c_\parallel/h_l \equiv u^*-c_\parallel/h_l^* 
\\
\xi_{-1} & = u_l^*-c_\perp/h_l^\star \equiv u_l^\sharp-c_\perp/h_l^\sharp
\\
\xi_0 & =  u^* \equiv \xi_0
\\
\xi_{+1} & = u_r^*+c_\perp/h_r^* \equiv u_r^\sharp+c_\perp/h_r^\sharp
\\
\xi_{+2} & = u_r+c_\parallel/h_r \equiv u^*+c_\parallel/h_r^*
\end{align*}
which are obviously compatible with the weak Riemann invariants of each wave.

\medskip

It remains to be seen 
whether the Riemann solver is actually ``entropy-consistent'', i.e. condition \eqref{solversufficientcondition} is satisfied.
To that aim, let us add two unknowns $\hat e_\parallel$, $\hat e_\perp$ to \eqref{eq:system1Deulerrelaxed}
such that
\begin{align}
\label{eq:epara}
\partial_{t}\left(h(\ux^2/2+\hat e_\parallel)\right) + \partial_{x}\left(h\ux(\ux^2/2+\hat e_\parallel^2)+\pi_\parallel\ux\right)
& = 0
\\
\label{eq:eperp}
\partial_{t}\left(h(\uy^2/2+\hat e_\perp)\right) + \partial_{x}\left(h\ux(\uy^2/2+\hat e_\perp)+\pi_\perp\uy\right) 
& = 0
\end{align}
hold.
On recalling the structure \eqref{eq:relaxationsolution} of solutions to Riemann problems,
if we use $(\hat e_\parallel)_o=e_\parallel(q_o)$ and $(\hat e_\perp)_o=e_\perp(q_o)$ for $o=l/r$ as left/right initial conditions,
then the following holds:
\begin{proposition}\label{prop:simpleconsistencycondition}
If the six following inequalities hold, for $o=l/r$,
\begin{align}
\label{cond1}
e_\parallel\left(q_o^*\right) 
& \le
\left(\hat e_\parallel\right)_o^*
\\
\label{cond2}
e_\perp\left(q_o^*\right) 
& \le
\left(\hat e_\perp\right)_o^*
\\
\label{cond3}
e_\perp\left(q_o^\sharp\right) 
& \le
\left(\hat e_\perp\right)_o^\sharp
\end{align}
then the 
entropy-consistency condition \eqref{solversufficientcondition} is satisfied\footnote{
  Note that it is of course not strictly necessary that \eqref{cond1} and \eqref{cond2} hold independently from one another,
  in fact it is sufficient that
  \begin{equation}
  \label{cond12}
  \left(e_\parallel + e_\perp\right)\left(q_o^*\right) 
  \le
  \left(\hat e_\parallel + \hat e_\perp\right)_o^*
  \end{equation}
  holds but it is easier to check \eqref{cond1} and \eqref{cond2} separately.
} with flux $\tilde G = [h\ux(\ux^2/2+\hat e_\parallel^2+\uy^2/2+\hat e_\perp)+\pi_\parallel\ux+\pi_\perp\uy]_{x/t=0}$ (a $c_0$ Riemann invariant).
\end{proposition}

\subsection{Choosing entropy-consistent relaxation parameters}
\label{sec:parameters}

We now explain how to satisfy the entropy-consistency conditions of Prop.~\ref{prop:simpleconsistencycondition}
for the relaxation parameters $c_\parallel^2,c_\perp^2,a^2$.

Condition \eqref{cond1} is classically satisfied provided the following 
condition
\begin{equation}
 \label{cond1bis}
 h^2\partial_h|_{h\sx^{-1/2},h\sz^{+1/2}} P_\parallel = h^2\left( gh + G h(3\sx+\sz) \right) \le c_\parallel^2
\end{equation}
is satisfied for all $h$ in between $h_o$ and $h_o^*$, $o=l/r$, see \cite{bouchut-2004,bouchut-boyaval-2013}.
The first step to show \eqref{cond1bis} is to compute $(\hat e_\parallel)_o^*=(\hat e_\parallel)_o^\sharp$ from the equation \eqref{eq:epara},
which can be done from the system in Lagrange coordinates augmented by the equation
$$
\partial_t (u^2/2+\hat e_\parallel) + \partial_x (u \pi_\parallel) = 0
$$
on noting 
that the following equation (in Lagrange coordinates) holds
$$
\partial_t (u^2/2+\pi_\parallel^2/2c_\parallel^2) + \partial_x (u \pi_\parallel) = 0
$$
so $(\hat e_\parallel)_o^*=(\hat e_\parallel)_o^\sharp$ can be obtained from the solution of 
$$
\partial_t (\hat e_\parallel-\pi_\parallel^2/2c_\parallel^2) = 0 \,.
$$
The second step subtracts $(P_\parallel\left(q_o^*\right)-\pi_\parallel^*)^2/(2c_\parallel^2)\ge0$
in the RHS of \eqref{cond1} rewritten
\beq
\label{cond1a}
\left(e_\parallel - \frac{P_\parallel^2}{2c_\parallel^2}\right)\left(q_o^*\right) 
\le
\left(e_\parallel - \frac{P_\parallel^2}{2c_\parallel^2}\right)\left(q_o\right) 
+ \frac{(\pi_\parallel^*)^2-P_\parallel\left(q_o^*\right)^2}{2c_\parallel^2}
\eeq
and uses the Riemann invariant $h^{-1}+\pi_\parallel/c_\parallel^2$ to show that, in fact,
\begin{equation*}
\label{cond1b}
\left(e_\parallel - \frac{P_\parallel^2}{2c_\parallel^2}\right)\left(h_o^*\right) 
\le
\left(e_\parallel - \frac{P_\parallel^2}{2c_\parallel^2}\right)\left(h_o\right) 
-P_\parallel(h^*_o)\left(\frac1{h^*_o}-\frac1{h_o}+\frac{P_\parallel(h^*_o)-P_\parallel(h_o)}{c_\parallel^2}\right)  
\end{equation*}
is enough, and therefore \eqref{cond1bis} after looking at variations in $h_o^*$.

\medskip

Similarly, the Riemann invariants $\hat e_\perp-\pi_\perp^2/2c_\parallel^2$ and $\hat e_\perp-\pi_\perp^2/2c_\perp^2$
give $(\hat e_\perp)_o^*$ and $(\hat e_\perp)_o^\sharp$.
Then, on noting \eqref{eq:sysstar}, a sufficient condition for \eqref{cond2} to hold reads 
\begin{multline}
\label{cond2b}
\left(e_\perp - \frac{P_\perp^2}{2c_\parallel^2}\right)\left(q_o^\star\right)
\le
\left(e_\perp - \frac{P_\perp^2}{2c_\parallel^2}\right)\left(q_o\right)
\\
-P_\perp(q^*_o)\left( 
 \frac{P_\perp(q^*_o)-P_\perp(q_o)}{c_\parallel^2}
 - \left(\frac{\sxy}{a^2h}\right)_o^\star + \left(\frac{\sxy}{a^2h}\right)_o
\right)  
\end{multline}
which rewrites with 
$\tilde G_o:=\left(\frac{G\sx}{h^2c_\parallel^2}\right)_o=\left(\frac{G\sx}{h^2c_\parallel^2}\right)_o^*$
and $\tilde \alpha_o:=\left(\frac{\sx}{h^2a^2}\right)_o=\left(\frac{\sx}{h^2a^2}\right)_o^*$
\begin{multline}
 \label{cond2bis}
 \left( 1 - \tilde G_o (h_o^\star)^4 \right) \left( \left(\frac{\sxy}h\right)_o^\star \right)^2
 \le 
 \left( 1 - \tilde G_o h_o^4 \right) \left( \left(\frac{\sxy}h\right)_o \right)^2 
\\
 - 2 (h_o^\star)^2 \left( \frac{\sxy}h \right)_o^*
 \Bigg( (\tilde G_o (h_o^\star)^2 - \tilde \alpha_o ) \left( \frac{\sxy}h \right)_o^* 
  - ( \tilde G_o h_o^2 - \tilde \alpha_o ) \left(\frac{\sxy}h\right)_o \Bigg)
\end{multline}
for $o=l/r$.
Now, the case $b=0$ is trivially satisfied since $({\sxy}/h)_o^*=({\sxy}/h)_o$.
Otherwise, when $b\neq0$, for $c_\parallel$ given the closed $q_1$ subsystem can be solved so that $h_o^\star$ is also fixed,
and \eqref{cond2bis} amounts to controlling the sign of a quadratic polynomial function of $({\sxy}/h)_o^*$
through $a$ and $c_\perp$.
So, if we ensure 
$$
1 + \tilde G_o (h_o^\star)^4 - 2 \tilde \alpha_o (h_o^\star)^2 < 0
\Leftrightarrow
\frac1{(h_o^\star)^2 } + \tilde G_o (h_o^\star)^2 < 2 \tilde \alpha_o \equiv 2 \left(\frac{\sx}{h^2}\right)_o \frac1{a^2}
$$
with $a$ small enough and if we next choose $({\sxy}/h)_o^*$ large enough (in magnitude) simply by controlling $c_\perp$
(it suffices to choose $c_\parallel^2-c_\perp^2>0$ small enough, given $c_\parallel$, $b$ and $a$) then \eqref{cond2bis} holds, thus \eqref{cond2}.

\medskip

Condition \eqref{cond3} can also be analyzed similarly to \eqref{cond1} 
and \eqref{cond2} on noting \eqref{eq:syssharp}.
Recalling the Riemann invariants for $c_{-1}$ and $c_1$, it suffices that for $o=l/r$
\begin{multline}
\label{cond3b}
\left(e_\perp - \frac{P_\perp^2}{2c_\perp^2}\right)\left(q_o^\sharp\right)
-
\left(\hat e_\perp - \frac{\pi_\perp^2}{2c_\perp^2}\right)_o^\star
\\
+ P_\perp(q^\sharp_o)\left( 
 \frac{P_\perp(q^\sharp_o)-(\pi_\perp)_o^*}{c_\perp^2}
 - \left(\frac{\sxy}{a^2h}\right)_o^\sharp + \left(\frac{\sxy}{a^2h}\right)_o^\star
\right) \le 0
\end{multline} 
holds. 
With $e_\perp(\sxy/h)=G(h^2/\sx)(\sxy/h)^2/2$ and $P_\perp(\sxy/h)=Gh^2(\sxy/h)$, note that
the LHS in \eqref{cond3b} is a quadratic polynomial in $(\frac{\sxy}h)^\sharp$
\begin{multline}
 \label{cond3bis}
 \left( 1 + \hat G_o - 2\sx/a^2 \right) \left( \left(\frac{\sxy}h\right)_o^\sharp \right)^2
- \left(\frac{2\sx}{Gh^2}\right)_o^\star \left(\hat e_\perp - \frac{\pi_\perp^2}{2c_\perp^2}\right)_o^\star
\\
 + 2 \sx \left( \frac1{a^2} \frac{\sxy}{h} - \frac{\pi_\perp}{c_\perp^2} \right)_o^* \left(\frac{\sxy}h\right)_o^\sharp
  \le 0
\end{multline}
where $\hat G_o:=\left(\frac{Gh^2\sx}{c_\perp^2}\right)_o^*$.
Now, to ensure condition \eqref{cond3b} similarly to \eqref{cond2b} i.e.
$$
\left( 1 + \hat G_o - 2\sx/a^2 \right)_o^* < 0
\Leftrightarrow 
c_\perp^2 + G(\sx h^2)_o^* \le \frac{2(\sx)_o^* c_\perp^2}{a^2}
$$
one cannot anymore choose $a$ and $c_\parallel^2-c_\perp^2>0$ independently small.
But \eqref{cond3b} can still be satisfied 
under the more stringent condition
$$ 
\frac{a^2}{c_\perp^2} \le \frac{2(\sx)_o^*}{c_\parallel^2 + G(h^2\sx)_o^*} \le \frac{2(\sx)_o^*}{c_\perp^2 + G(h^2\sx)_o^*} 
$$
plus simultaneously a large enough $({\sxy}/h)_o^\sharp$ (in magnitude).

If $b(u^*-u_o)\neq0$ then one can again ensure a large enough $({\sxy}/h)_o^\sharp$ with $c_\parallel^2-c_\perp^2>0$ small,
which is obviously compatible with our requirements for condition \eqref{cond2b}.

If $b(u^*-u_o)=0$, then \eqref{cond3} simplifies 
and reads
\begin{multline}
\label{eq:bduzero}
\frac{G}2\left(\frac{h^2}{\sx}\right)^*_o \left( \left|\frac{\sxy}h\right|_{\sharp,o}^2 - \left|\frac{\sxy}h\right|_o^2 \right) 
\le \frac{G(h^*_o)^2}{a^2} \left(\frac{\sxy}h\right)_{\sharp,o} \left( \left(\frac{\sxy}h\right)_{\sharp,o} - \left(\frac{\sxy}h\right)_o \right)
\\
- \frac{\left(G(h^*_o)^2\left(\frac{\sxy}h\right)_{\sharp,o}-Gh_o^2\left(\frac{\sxy}h\right)_o\right)^2}{2c_\perp^2}
\end{multline}
which is satisfied if $c_\perp^2$ is large enough.
Now, recalling the constraint above on $c_\perp^2/c_\parallel^2<1$, this requires one to increase $c_\parallel$.
So an entropy-consistent choice for the relaxation parameters can always be identified by a simple algorithm explained in the next Section~\ref{sec:scheme}
with the whole computational scheme.
\mycomment{But the latter may be a problem for CFL, except if this is trivial.}


In the end, we note that 
the parameter $b=2G\ro^2\sxy$, which has no sign a priori,
can be choosen ``freely''.
To minimize numerical diffusion, we compute it as the solution to $\partial_t b + u \partial_x b = 0$
initialized consistently with its definition. 

\subsection{Numerical schemes for SVTM and SVUCM}
\label{sec:scheme}

Finally, let us summarize our numerical scheme to simulate SVTM, which turns out to be useful for SVUCM at the cost of very few modifications.

\smallskip

But first, for each Riemann problem, we propose to also ``relax'' $a$, $c_\parallel$ and $c_\perp$ as independent variables solutions to pure transport equations,
similarly to $b$, for the sake of more precision. 
Then, the exact values of the intermediate states change, and become a function of \emph{left and right} parameter values.
For SVTM 1D Riemann problems, they become:
\begin{equation}
\label{final1}
{\tiny 
\hspace{-2cm}
\begin{array}{rl}
u^* = \frac{c_{\parallel,l}u_l+\pi_{\parallel,l}+c_{\parallel,r}u_r-\pi_{\parallel,r}}{c_{\parallel,l}+c_{\parallel,r}},
&\qquad
\pi^*_\parallel = \frac{c_{\parallel,r}(\pi_{\parallel,l}+c_{\parallel,l}u_l)+c_{\parallel,l}(\pi_{\parallel,r}-c_{\parallel,r}u_r)}{c_{\parallel,l}+c_{\parallel,r}},
\\
\frac{1}{h_{\parallel,l}^*}=\frac{1}{h_{\parallel,l}}+\frac{c_{\parallel,r}(u_r-u_l)+\pi_{\parallel,l}-\pi_{\parallel,r}}{c_{\parallel,l}(c_{\parallel,l}+c_{\parallel,r})},
&\quad
\frac{1}{h_{\parallel,r}^*}=\frac{1}{h_{\parallel,r}}+\frac{c_{\parallel,l}(u_r-u_l)+\pi_{\parallel,r}-\pi_{\parallel,l}}{c_{\parallel,r}(c_{\parallel,l}+c_{\parallel,r})},
\\
v_l^\star = v_l+\frac{b_l}{c_{\perp,l}^2-c_{\parallel,l}^2}
\frac{c_{\parallel,r}(u_l-u_r)+\pi_{\parallel,r}-\pi_{\parallel,l}}{c_{\parallel,l}+c_{\parallel,r}},
&\quad
v_r^\star = v_r+\frac{b_r}{c_{\perp,r}^2-c_{\parallel,r}^2}
\frac{c_{\parallel,l}(u_r-u_l)+\pi_{\parallel,r}-\pi_{\parallel,l}}{c_{\parallel,l}+c_{\parallel,r}},
\\
\pi_{\perp,l}^\star 
= \pi_{\perp,l} - \frac{b_l c_{\parallel,l}}{c_{\perp,l}^2-c_{\parallel,l}^2}
\frac{c_{\parallel,r}(u_l-u_r)+\pi_{\parallel,r}-\pi_{\parallel,l}}{c_{\parallel,l}+c_{\parallel,r}},
&\quad
\pi_{\perp,r}^\star 
= \pi_{\perp,r} + \frac{b_r c_{\parallel,r}}{c_{\perp,r}^2-c_{\parallel,r}^2}
\frac{c_{\parallel,l}(u_r-u_l)+\pi_{\parallel,r}-\pi_{\parallel,l}}{c_{\parallel,l}+c_{\parallel,r}},
\\
\left(\frac{\sxy}h\right)_l^\star 
= \left(\frac{\sxy}h\right)_l - \frac{a^2_l b_l / c_{\parallel,l}}{c_{\perp,l}^2-c_{\parallel,l}^2}
\frac{c_{\parallel,r}(u_l-u_r)+\pi_{\parallel,r}-\pi_{\parallel,l}}{c_{\parallel,l}+c_{\parallel,r}},
&\quad
\left(\frac{\sxy}h\right)_r^\star 
= \left(\frac{\sxy}h\right)_r + \frac{a^2_r b_r / c_{\parallel,r}}{c_{\perp,r}^2-c_{\parallel,r}^2}
\frac{c_{\parallel,l}(u_r-u_l)+\pi_{\parallel,r}-\pi_{\parallel,l}}{c_{\parallel,l}+c_{\parallel,r}},
\\
v^\sharp = \frac{c_{\perp,l}v_l^*+\pi_{\perp,l}^*+c_{\perp,r}v_r^*-\pi_{\perp,r}^*}{c_{\perp,l}+c_{\perp,r}}
&\qquad
\pi_{\perp}^\sharp = \frac{c_{\perp,r}(\pi_{\perp,l}^*+c_{\perp,l}v_l^*)+c_{\perp,l}(\pi_{\perp,r}^*-c_{\perp,r}v_r^*)}{c_{\perp,l}+c_{\perp,r}}
\\
\left(\frac{\sxy}h\right)_r^\sharp 
= \left(\frac{\sxy}h\right)_r^\star + \frac{a^2_r}{c_{\perp,r}^2}
\frac{c_{\perp,r}c_{\perp,l}(v_l^*-v_r^*)+c_{\perp,l}(\pi_{\perp,r}^*-\pi_{\perp,l}^*)}{c_{\perp,l}+c_{\perp,r}}
&\quad
\left(\frac{\sxy}h\right)_l^\sharp 
= \left(\frac{\sxy}h\right)_l^\star + \frac{a^2_l}{c_{\perp,l}^2}
\frac{c_{\perp,r}c_{\perp,l}(v_l^*-v_r^*)+c_{\perp,r}(\pi_{\perp,l}^*-\pi_{\perp,r}^*)}{c_{\perp,l}+c_{\perp,r}}
\\
(\ro^{-2}\sx)_l^\star = (\ro^{-2}\sx)_l^\sharp = \ro^{-2}\sx
&\quad
(\ro^{-2}\sx)_r^\star = (\ro^{-2}\sx)_r^\sharp = \ro^{-2}\sx
\\
(\ro^{2}\sz)_l^\star = (\ro^{2}\sz)_l^\sharp = \ro^{2}\sz
&\quad
(\ro^{2}\sz)_r^\star = (\ro^{2}\sz)_r^\sharp = \ro^{2}\sz
\\
\end{array}
}
\end{equation}
however the 
conditions (\ref{cond1}--\ref{cond2}--\ref{cond3}) do not change,
and our discussion in Section~\ref{sec:parameters} 
to achieve discrete entropy dissipation by well-chosen $a$, $c_\parallel$ and $c_\perp$ is still valid
(although the identification of numerical values for the right and left relaxation parameters may become more difficult).

\smallskip

Second, to use a similar approach for SVUCM, note also that \eqref{eq:system1D2} becomes
\beq
\label{eq:system1D2uc}
\left\lbrace
\begin{aligned}
\partial_{t}    \ro 
+ \partial_{x}( \ro \ux ) 
& 
= 0
\\
\partial_{t}(  \ro\ux  )
+ \partial_{x}( \ro \ux^2 ) 
+ \partial_x( g \ro^2/2 ) 
+ \partial_x\left( G \ro(\sz -\sx) \right)
&
= 0
\\
\partial_{t}( \ro^{2}\sx )  + \ux \partial_{x}( \ro^{2}\sx ) 
&
= 0
\\
\partial_{t}(  \ro^{-2}\sz  ) + \ux\partial_{x}( \ro^{-2}\sz  ) 
&
= 0
\\
\partial_{t}( \sy-\sxy^2/\sx  ) + \ux\partial_{x}( \sy-\sxy^2/\sx  ) 
&
= 0
\\
\partial_{t}( \ro \uy )
+ \partial_{x}( \ro \ux \uy ) 
- \partial_x\left( G \ro\sxy \right)
&
= 0
\\
\partial_{t}( \ro\sxy ) + \ux\partial_{x}( \ro\sxy ) - \ro\sx \partial_x\uy
&
= 0
\end{aligned}
\right.
\eeq
while
the pressures associated with the energy contributions $E_\parallel$ and $E_\perp$, which are exactly the same as in SVTM,
then respectively read $P_\parallel=\frac{gh^2}2 + Gh(\sz-\sx)$ and $P_\perp=-G\ro\sxy$.
Now, the pressures follow similar evolution equations for smooth solutions
with parameters $c_\parallel^2=h^2\left(gh + G(3\sz+\sx)\right)$, $c_\perp^2=Gh^2\sx$, $b\equiv0$,
and this is why we can use the same (pseudo-)relaxation approach.
For SVUCM, we define as Riemann solver the (exact) solution to the following linearly degenerate system 
\beq
\label{eq:system1Duceulerrelaxed}
\left\lbrace
\begin{aligned}
\partial_{t} \ro - \partial_{x} (\ro\ux)
& 
= 0
\\
\partial_{t} (\ro\ux) + \partial_x (\ro\ux^2+\pi_\parallel)
&
= 0
\\
\partial_t (\ro\pi_\parallel) + \partial_x (\ro\ux\pi_\parallel) + c_\parallel^2 \partial_x \ux
&
= 0
\\
\partial_{t} (-\ro^2\sxy) + \partial_x (-\ro^2\ux\sxy) + a^2 \partial_x \uy
&
= 0
\\
\partial_{t}(\ro\uy) + \partial_x (\ro\ux\uy+\pi_\perp)
&
= 0
\\
\partial_t (\ro\pi_\perp) + \partial_x (\ro\ux\pi_\perp) 
+ c_\perp^2 \partial_x \uy
&
= 0
\\
\partial_{t}( \ro^{-1}\sz ) + \partial_{x}( \ro^{-1}\sz\ux )  
&
= 0
\\
\partial_{t}(  \ro^3\sx  ) + \partial_{t}(  \ro^3\sx\ux  )
&
= 0
\\
\partial_{t}\left(  \ro( \sy-\sxy^2/\sx  )  \right) + \partial_{x}\left(  \ro( \sy-\sxy^2/\sx  )\ux  \right)
&
= 0
\end{aligned}
\right.
\eeq
which is easily deduced from \eqref{eq:system1Duceulerrelaxed} on noting the new value of $a^2=\sx\ro^2$ that is consistent with its definition in smooth cases.

As opposed to SVTM, it always holds $b=0$ in the SVUCM case.
This is consistent with the fact that there is no non-conservative product between nonlinear fields in \eqref{eq:system1D2uc} as opposed to \eqref{eq:system1D2}
(observe indeed that the subsystem for $q_1$ is closed in SVUCM, like ine SVTM, and that for $q_2=(\sxy,v,\sy)$ is also closed -- unlike SVTM).
Then, as a consequence, entropy-consistency can be ensured for SVUCM following the same approach as for SVTM !

In fact, since $b=0$, only the two conditions \eqref{cond1} and \eqref{cond3} have to be checked.
The first one is still satisfied with the choice in \cite{bouchut-boyaval-2013}.
The second 
reads
\begin{multline}
\label{eq:bzero}
\frac{G}2\left(\frac1{h^2\sx}\right)^*_o \left( \left|h \sxy\right|_{\sharp,o}^2 - \left|h \sxy\right|_o^2 \right) 
\le \frac{G}{a^2} \left(h\sxy\right)_{\sharp,o} \left( \left(h\sxy\right)_{\sharp,o} - \left(h\sxy\right)_o \right)
\\
- \frac{G}{2c_\perp^2}\left(\left(h\sxy\right)_{\sharp,o}-\left(h\sxy\right)_o\right)^2
\end{multline}
since $b=0$ and is always satisfied when $a^2\le\sx\ro^2$ 
(unlike \eqref{eq:bduzero} in SVTM).

\smallskip

Now, given a polygonal mesh of $\R^2$, with cells $V_i$ and interfaces $\Gamma_{ij}$,
time evolutions of viscoelastic flows can be numerically simulated 
using a piecewise-constant state vector
$$
 q:=(H,HU,HV,HC_{xx},HC_{yy},HC_{xy}/\sqrt{C_{xx}C_{yy}},HC_{zz})
$$
with convex admissible domain
$$
\{H,HC_{xx},HC_{yy},HC_{zz}>0 \;;\quad -H<HC_{xy}/\sqrt{C_{xx}C_{yy}}<H\} \,,
$$
that is a Finite-Volume (FV) approximation either for SVTM or for SVUCM. 

Time-discrete sequences can be computed in two sub-steps.
First, define $q_i^{n+1,-}$ in \eqref{eq:conservation}
with a numerical flux of the form (\ref{eq:rotationinvariantflux}--\ref{eq:numfluxHLL})
using the 1D Riemann solver $R$ 
defined by \eqref{final1} for SVTM, and the straightforward modification proposed above for SVUCM.
Second, compute \eqref{cnvxcombinationu}, \eqref{cnvxcombinationh} and \eqref{cnvxcombinationz}.


\begin{proposition}
At each time step, 
there exist left and right initial states of the relaxation parameters $c_\parallel$, $c_\perp$ and $b$ for all Riemann problems at interfaces $\Gamma_{ij}$
such that 
the convex admissible domain for $q$ is preserved,
and the second-principle formulation \eqref{eq:discreteentropy2} is satisfied under the stringent CFL condition \eqref{eq:CFL}.
\end{proposition}


In practice, we need \emph{numerical values} of relaxation parameters 
such that conditions \eqref{cond1} and \eqref{cond3} hold, plus \eqref{cond2} for SVTM: they are useful in \eqref{final1}.
An adequate choice of the relaxation parameters can be computed in each Riemann problem as follows, for instance.

First, initialize $c_\parallel$ to ensure 
\eqref{cond1} 
like \cite{bouchut-boyaval-2013} for 1D cases
{\small
\begin{align}
\label{eq:initialization1}
c_{\parallel,l} = h_{\parallel,l}\left(\sqrt{\partial_hP_\parallel(q_{l})}+2\left((u_l-u_r)_++\frac{(\pi_{\parallel,r}-\pi_{\parallel,l})_+}{h_{\parallel,l}\sqrt{\partial_hP_\parallel(q_{l})}+h_{\parallel,r}\sqrt{\partial_hP_\parallel(q_{r})}}\right)\right) \,,
\\
\label{eq:initialization2}
c_{\parallel,r} = h_{\parallel,r}\left(\sqrt{\partial_hP_\parallel(q_{r})}+2\left((u_l-u_r)_++\frac{(\pi_{\parallel,l}-\pi_{\parallel,r})_+}{h_{\parallel,l}\sqrt{\partial_hP_\parallel(q_{l})}+h_{\parallel,r}\sqrt{\partial_hP_\parallel(q_{r})}}\right)\right) \,.
\end{align}
}
where we denoted
$ 
\partial_hP_\parallel(q_o) = \partial_h|_{h\sx^{-1/2},h\sz^{+1/2}}P_\parallel(q_o) = g h_o + G(3\sx+\sz)_o \,.
$

Given $c_{\parallel,o}$, only $a^2_o = c_{xx,o}h_o^2$ remains to compute for SVUCM, while
$$
a^2_o = c_{xx,o}\wedge r^0_o \frac{2c_{xx,o}^*}{1+Gh_o^2c_{xx,o}/c_{\parallel,o}^2}
$$
for SVTM is not enough.

For SVTM, one still needs to compute $c_{\perp,o}^2=r^n_o c_{\parallel,o}^2$ 
after $n$ iterations of a strictly increasing sequence $r^n_o<r^{n+1}_o\le 1$ with limit $1$
starting from 
$r^0_o=Gh_o^2c_{xx,o}/c_{\parallel,o}^2>0$ 
until \eqref{cond2} (or the weaker conditions \eqref{cond1}+\eqref{cond2}) 
are satisfied for some $n$ \emph{both} for $o=l/r$ simultaneously. 
Moreover, given $c_{\parallel,o}^2,a^2_o$, $b_o=2Gh_o^2c_{xy,o}$ and 
$c_{\perp,o}^2/c_{\parallel,o}^2=r^n_o\le 1$, one still has to inspect \eqref{cond3} in SVTM case.
If $b_o(u_o-u^*)\neq0$ we iterate further on the sequence $r^n_o<r^{n+1}_o\le 1$ defining $c_{\perp,o}^2$
until \eqref{cond3} are satisfied for some $n$ \emph{both} for $o=l/r$ simultaneously.
Otherwise, if $b_o(u_o-u^*)=0$ for any $o=l/r$, then $c_{\parallel,o}$ has to be increased where \eqref{cond3} is not satisfied:
one can use some unbounded increasing sequence $\tilde r_n$ and iterate on $n$ to define $c_{\parallel,o}$ starting with (\ref{eq:initialization1}--\ref{eq:initialization2}).

\smallskip

Finally, using the analytical Riemann solution of section \ref{sec:oned} one has
\begin{multline}
\label{finalflux}
\frac1{\tau^n}\int_{-\left(\sum_j \frac{|\Gamma_{ij}|}{|V_i|}\right)^{-1}}^0 R(x/\tau^n,q_l,q_r) dx
=  \left( \xi_{-2,-} + \bar s_i^n \right) q_l 
+  \left( \xi_{-1,-} - \xi_{-2,-} \right) q_l^* 
\\
+  \left( \xi_{0,-} - \xi_{-1,-} \right) q_l^\sharp 
+  \left( \xi_{1,-} - \xi_{0,-} \right) q_r^\sharp 
+  \left( \xi_{2,-} - \xi_{1,-} \right) q_r^*
+  \left( - \xi_{2,-} \right) q_r
\\
=  \bar s_i^n q_l 
+  \left( \xi_{-1,-} - \xi_{-2,-} \right) (q_l^*-q_l) 
+  \left( \xi_{0,-} - \xi_{-1,-} \right) (q_l^\sharp-q_l) 
\\
+  \left( \xi_{1,-} - \xi_{0,-} \right) (q_r^\sharp-q_l) 
+  \left( \xi_{2,-} - \xi_{1,-} \right) (q_r^*-q_l)
+  \left( - \xi_{2,-} \right) (q_r-q_l)
\end{multline}
where we have denoted $x_-=\min(x,0)$ the non-positive part of a real number $x$
and $\bar s_i^n:=\frac1{\tau^n}\left(\sum_j \frac{|\Gamma_{ij}|}{|V_i|}\right)^{-1}\ge s_i^n>0$.
%
%
Note that \eqref{finalflux} is 
consistent with conservation laws for ``conservative'' components of variable $q$ like $U$ i.e.
$$ \bF_{i\to j}^U(q_i,q_j;\bn_{i\to j}) + \bF_{j\to i}^U(q_j,q_i;\bn_{j\to i}) = 0 \,. $$ 
Therefore, under our stringent CFL condition \eqref{eq:CFL}, 
\eqref{eq:cvxcombin2}
rewrites 
\begin{multline}
\label{eq:cvxcombin3}
q_i^{n+1,-} = \sum_j \frac{|\Gamma_{ij}|}{|V_i|} \Oij \int_{-\left(\sum_j \frac{|\Gamma_{ij}|}{|V_i|}\right)^{-1}}^0 R(x/\tau^n,\Oij^{-1}q_i^n,\Oij^{-1}q_j^n) dx
\\
= q_i^n + \tau^n \sum_j \frac{|\Gamma_{ij}|}{|V_i|} 
\Bigg( \left( \xi_{-1,-} - \xi_{-2,-} \right) (\Oij q_l^*-q_i^n) 
\\
+  \left( \xi_{0,-} - \xi_{-1,-} \right) (\Oij q_l^\sharp-q_i^n) 
+  \left( \xi_{1,-} - \xi_{0,-} \right) (\Oij q_r^\sharp-q_i^n) 
\\
+  \left( \xi_{2,-} - \xi_{1,-} \right) (\Oij q_r^*-q_i^n)
+  \left( - \xi_{2,-} \right) (q_j^n-q_i^n)
\Bigg)
\end{multline}
where $q^{(*,\sharp)}_o$ are the intermediate states of the Riemann solver $R(\xi,\Oij^{-1}q_i^n,\Oij^{-1}q_j^n)$,
$\xi_{-2}<\xi_{-1}<\xi_{0}<\xi_{+1}<\xi_{+2}$ its wave speeds, 
and where $\tau^n$ can be computed at each time step after 
all Riemann problems at all faces to satisfy \eqref{eq:CFL}.
%

\section{Numerical illustration 
and discussion} 
\label{sec:num}


We now illustrate the SVTM and SVUCM models using numerical solutions computed with our FV schemes in benchmark test cases.

\subsection{Stoker test case}
\label{sec:stoker}

This is a well-known benchmark test case for the (time-dependent, inviscid) Saint-Venant shallow-water equations, 
which models 
an idealized \emph{dam-break} 
(i.e. the propagation under gravity of 
a shock wave in a finite-depth fluid initially at rest) \cite{stoker-1957}.
A solution for $t\in(0,.2)$ is computed in a square $(x,y)\in[0,1]^2$ starting from the initial condition
$$
(H,U,V,C_{xx},C_{yy},C_{xy},C_{zz}) = 
\begin{cases}
(3,0,0,1,1,0,1) & x+y<1
\\
(1,0,0,1,1,0,1) & x+y>1 
\end{cases}
$$
that consists in two equilibrium rest states on each side of the line $x+y=1$.
For the boundary conditions, we use the ``ghost cell'' method (see e.g. \cite{leveque-2002}) assuming translation invariance along $x+y$ isolines.

The main point of Stoker test case is usually 
to compare 
the speed of the shock front with 
observations
(the solution expected -- numerically at least -- is indeed a free-shear flow). 
%
Alternatively, assuming translation invariance along $x+y$ isolines, 
this is a 1D Riemann problem that uniquely determines $q_1$ (as well as $q_2$ for SVUCM).
So the test case can be used to accurately understand 
the new variable $\bC$ in out-of-equilibirium dynamics, phenomenologically at least,
as a function of the viscoelastic parameters $G,\lambda$.
\mycomment{Note there are still experimental observations of surface waves in viscoelastic fluids \cite{MONROY20174}.}

Moreover, for 
benchmarking purposes, we also aim at comparing the new 2D scheme with former 1D numerical results obtained along $x=y$
in our previous work \cite{bouchut-boyaval-2013} for the subsystem $q_1$ of SVUCM model,
at Froude number $g^{-1/2}=.3$ ($g=10$), 
elasticiy number $G=10$ 
and Weissenberg 
number $\lambda=1$. 
This of course assumes that our 2D scheme converges to a translation-invariant solution.

\begin{figure}[hbtp]
\centering
\reflectbox{\rotatebox[origin=c]{180}{\includegraphics[scale=.45]{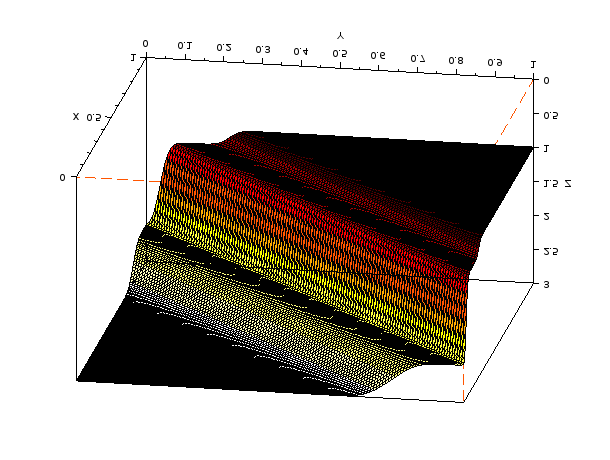}}}
\caption{\label{fig0} Stoker test case: flow depth $H(x,y$ 
at final time $T=.2$}
\end{figure}

\begin{figure}[hbtp]
\reflectbox{\rotatebox[origin=c]{180}{\includegraphics[scale=.4,clip,trim = 50 50 50 40]{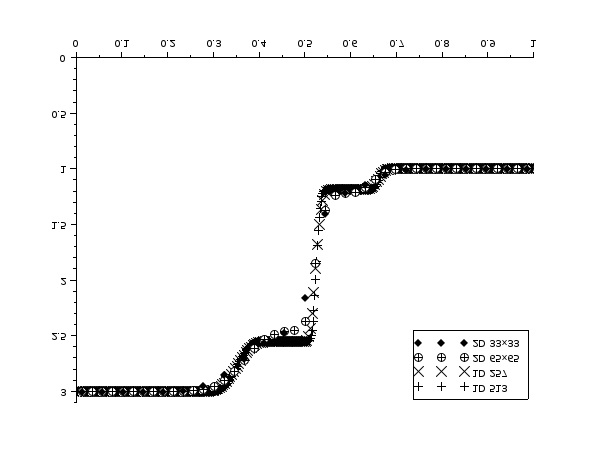}}}
\reflectbox{\rotatebox[origin=c]{180}{\includegraphics[scale=.4,clip,trim = 50 50 50 40]{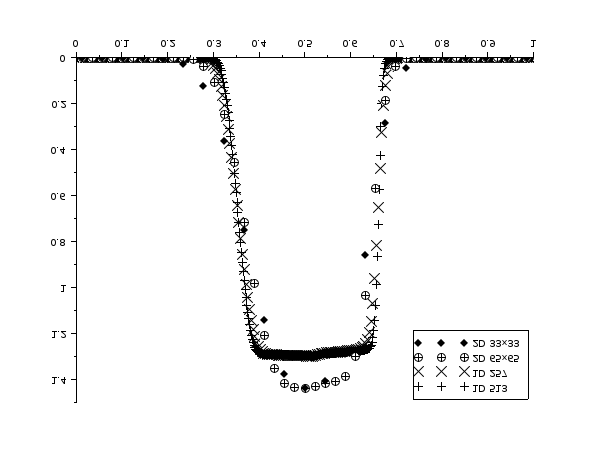}}}
\reflectbox{\rotatebox[origin=c]{180}{\includegraphics[scale=.4,clip,trim = 50 50 50 40]{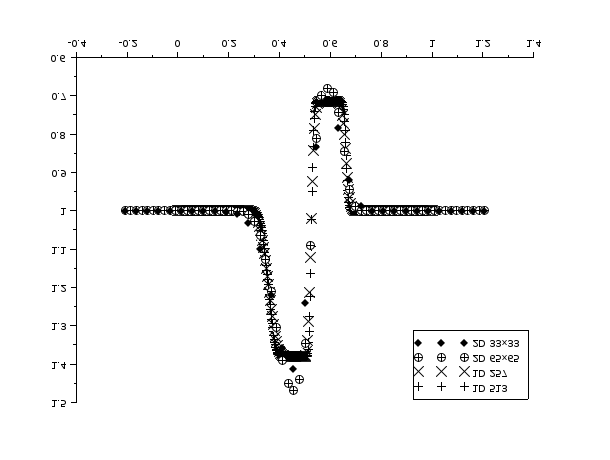}}}
\reflectbox{\rotatebox[origin=c]{180}{\includegraphics[scale=.4,clip,trim = 50 50 50 40]{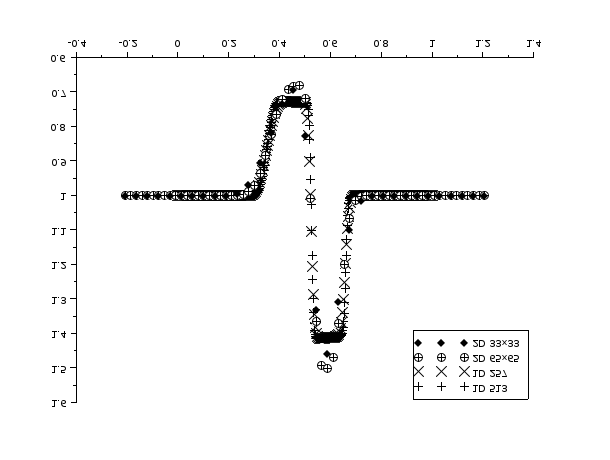}}}
\caption{\label{fig1} Stoker test case: $H,\bU\cdot\bn,C_{nn},C_{zz}$ computed by SVUCM 2D (cross-section $x=y$) 
and 1D (along $\bn$ normal to the initial discontinuity at $x+y=1$)}
\end{figure}

In Fig.~\ref{fig1}, this is exactly what we observe, see e.g. the variable $H$ in Fig.\ref{fig0}
(note that 
translation-invariance is used to define the values of ``ghost cells'' in Riemann problems at boundary faces).
We compare results obtained with 2D Cartesian meshes of 
$(2^5+1)=1089$, 
$(2^6+1)=4225$ 
cells and 1D regular grids with $2^8+1=257$ and $2^9+1=513$ cells.
The 2D solutions seem to preserve the initial translation invariance and converge to the (unique) 1D translation-invariant solution.
(The positions of the fronts seem already quite well resolved with our coarse mesh, at least,
although the 2D state values are only $10\%$ accurate relatively to 1D.)

\begin{figure}[hbtp]
\reflectbox{\rotatebox[origin=c]{180}{\includegraphics[scale=.4,clip,trim = 50 50 50 40]{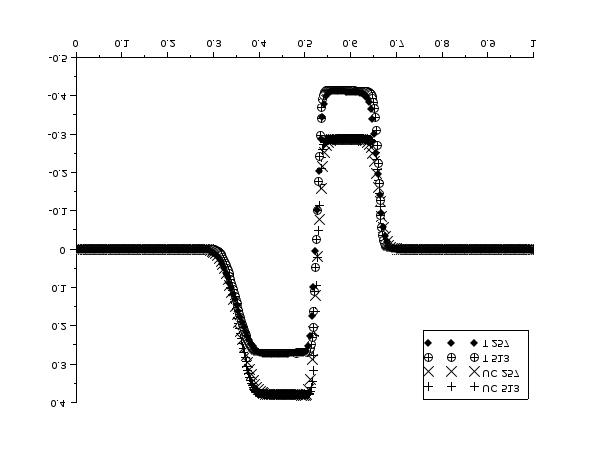}}}
\reflectbox{\rotatebox[origin=c]{180}{\includegraphics[scale=.4,clip,trim = 50 50 50 40]{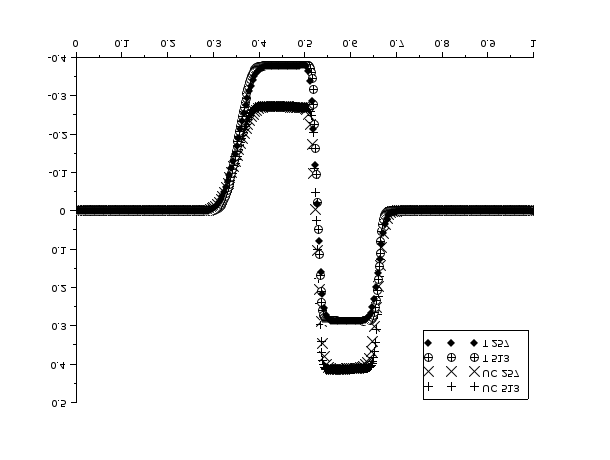}}}
\caption{\label{fig2} Stoker test case: $\Sigma_{nn},\Sigma_{zz}$ 
(left/right) for SVTM and SVUCM (T/UC) 
along the normal $\bn$ to the initial discontinuity at $x+y=1$ (in $G$ units).
}
\end{figure}

For 
benchmarking purposes, the test case also allows one to compare SVTM and SVUCM in a simple configuration.
We do not observe significant differences for $H,\bU$ at such short times, however
the non-zero stress components $\Sigma_{nn},\Sigma_{zz}$ are not the same (although they have similar tendencies)
see Fig.~\ref{fig2} for the various 1D (converged) values.

Last, variations in the parameters $G$ and $\lambda$ can also be well understood phenomenologically 
(from teh viscoelastic physics viewpoint) in the present 1D testcase.
Since analyzing variations in the parameters $G$ and $\lambda$ was already done in \cite{bouchut-boyaval-2013} for SVUCM,
we concentrate here on SVTM.

\begin{figure}[hbtp]
\reflectbox{\rotatebox[origin=c]{180}{\includegraphics[scale=.22,clip,trim = 50 50 50 40]{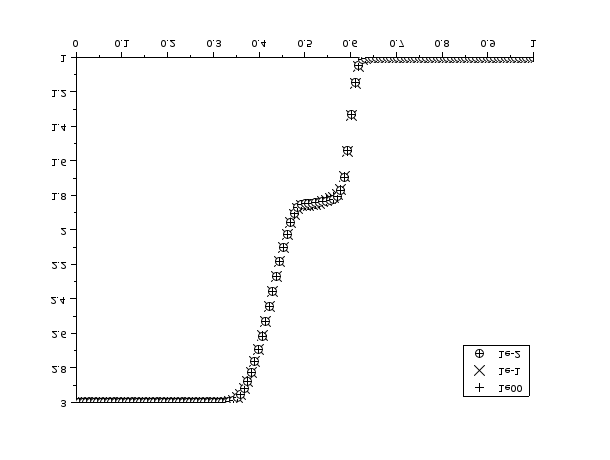}}}
\reflectbox{\rotatebox[origin=c]{180}{\includegraphics[scale=.22,clip,trim = 50 50 50 40]{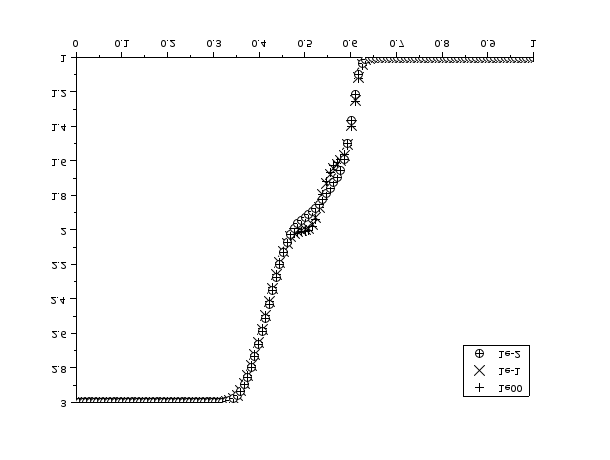}}}
\reflectbox{\rotatebox[origin=c]{180}{\includegraphics[scale=.22,clip,trim = 50 50 50 40]{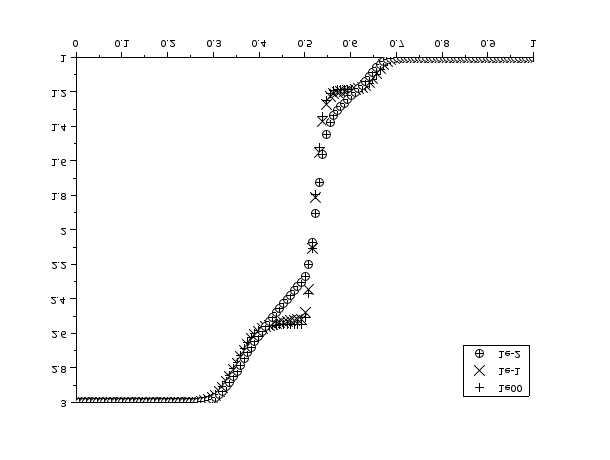}}}

\reflectbox{\rotatebox[origin=c]{180}{\includegraphics[scale=.22,clip,trim = 50 50 50 40]{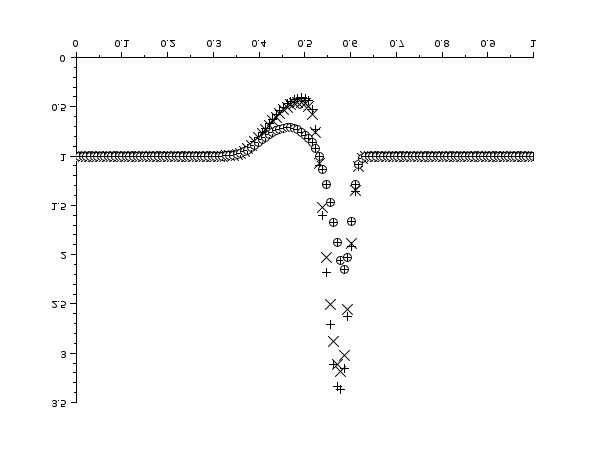}}}
\reflectbox{\rotatebox[origin=c]{180}{\includegraphics[scale=.22,clip,trim = 50 50 50 40]{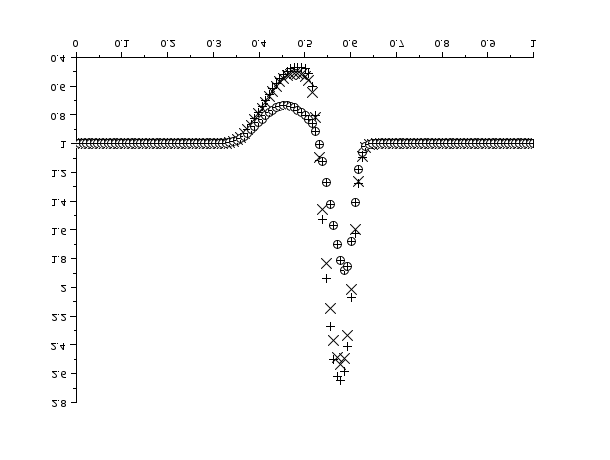}}}
\reflectbox{\rotatebox[origin=c]{180}{\includegraphics[scale=.22,clip,trim = 50 50 50 40]{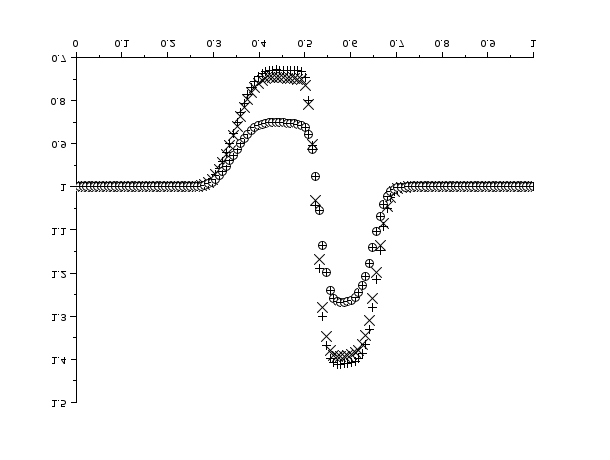}}}

\reflectbox{\rotatebox[origin=c]{180}{\includegraphics[scale=.22,clip,trim = 50 50 50 40]{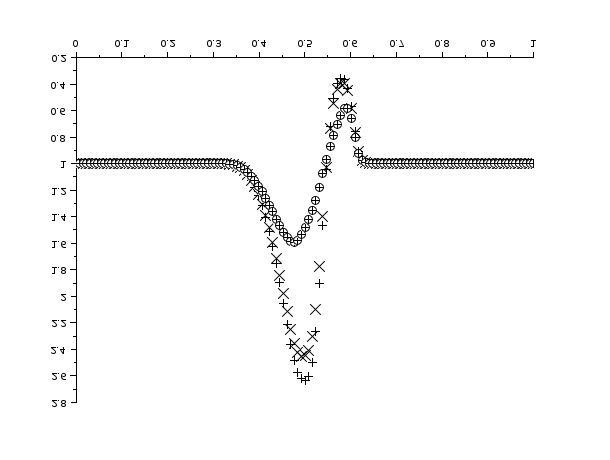}}}
\reflectbox{\rotatebox[origin=c]{180}{\includegraphics[scale=.22,clip,trim = 50 50 50 40]{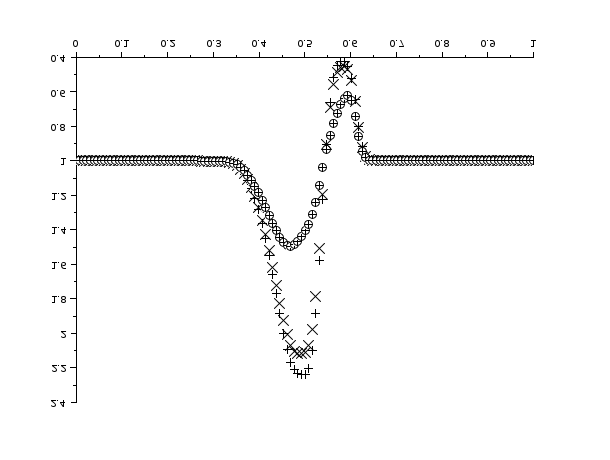}}}
\reflectbox{\rotatebox[origin=c]{180}{\includegraphics[scale=.22,clip,trim = 50 50 50 40]{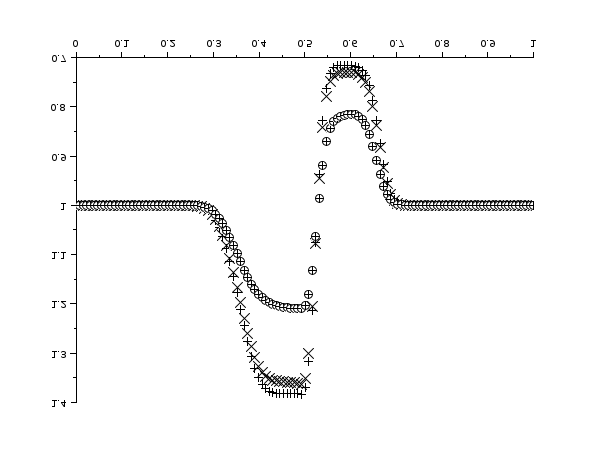}}}
\caption{\label{fig21} Stoker test case: $H$ (top), $C_{xx}$ (middle), $C_{zz}$ (bottom) at $T=.2$ for $\lambda=.01,.1,1$
when $G=.1,1,10$ (left, center and right) with $257$ cells.}
\end{figure}

In Fig.~\ref{fig21} we show $H,C_{xx},C_{zz}$ for SVTM at $T=.2$ when $G=.1,1,10$ and $\lambda=.01,.1,1$.
Increasing $G$ at fixed $\lambda$ only slightly increases the speeds of the genuinely nonlinear waves,
while it more essentially increases the jump in $H$ at the linearly degenerate wave (a contact discontinuity)
and decreases the jumps in $C_{xx},C_{zz}$.
This is physically coherent with the fact that the elasticity $G$ controls how difficult it is to locally deform 
the fluid materials of depth $H$, and connecting two equilibria at $H=3$ and $H=1$ through deformations becomes harder
as $G$ increases. 
However, if $\lambda$ simultaneously decreases, then variations in space of the strain $C_{xx},C_{zz}$ are fast smoothed back to equilibrium
and a viscous profil arises (see $G=10$, $\lambda=10^{-2}$ in Fig.~\ref{fig21}).
We recall that \emph{both} models SVUCM and SVTM formally converge to the viscous Saint-Venant equations in asymptotics $\lambda,G^{-1}\to0$
where $1/G\lambda$ remains bounded and defines a Reynolds number. 
Of course, at fixed parameter values, it is not so clear to define how close solutions are from a viscous approximation,
all the less when the Froude number $g^{-1/2}$ also varies.
This may be an interesting direction for future research directions, insofar as $g^{-1/2}$ should be quite low in real viscoelastic fluids.

\begin{figure}[hbtp]
\reflectbox{\rotatebox[origin=c]{180}{\includegraphics[scale=.3,clip,trim = 50 50 50 40]{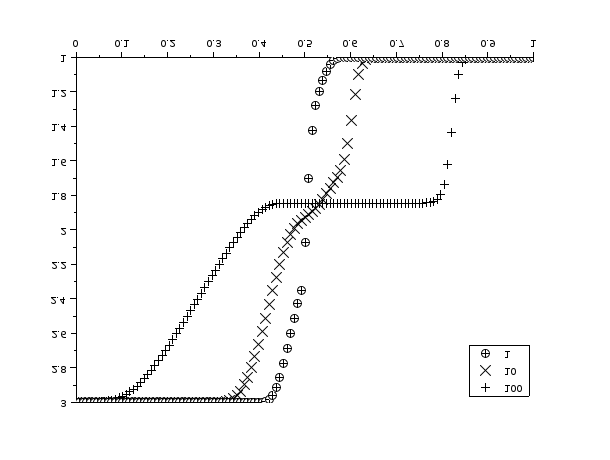}}}
\reflectbox{\rotatebox[origin=c]{180}{\includegraphics[scale=.3,clip,trim = 50 50 50 40]{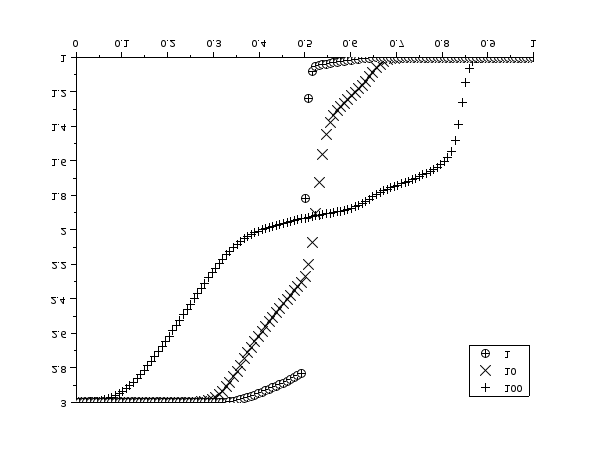}}}

\reflectbox{\rotatebox[origin=c]{180}{\includegraphics[scale=.3,clip,trim = 50 50 50 40]{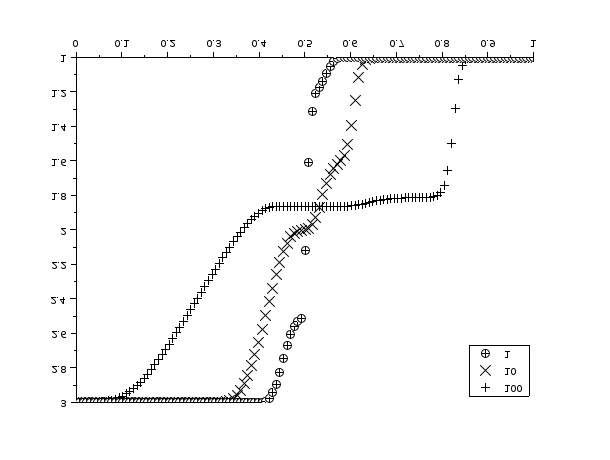}}}
\reflectbox{\rotatebox[origin=c]{180}{\includegraphics[scale=.3,clip,trim = 50 50 50 40]{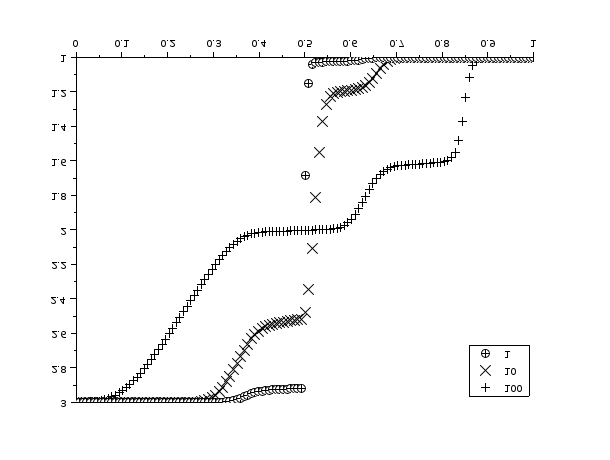}}}
\caption{\label{fig22} Stoker test case: $H$ at $T=.2$ for $g=1,10,100$ when $\lambda=.01,1$ (top/bottom) and $G=1,10$ (left/right) with $257$ cells.}
\end{figure}

\mycomment{TO DO: show tendency as a function of elasticity $G$ and Weissenberg $\lambda$}

\subsection{Collapse of a column}

\begin{figure}[hbtp]
\reflectbox{\rotatebox[origin=c]{180}{\includegraphics[scale=.27,clip,trim = 70 50 50 40]{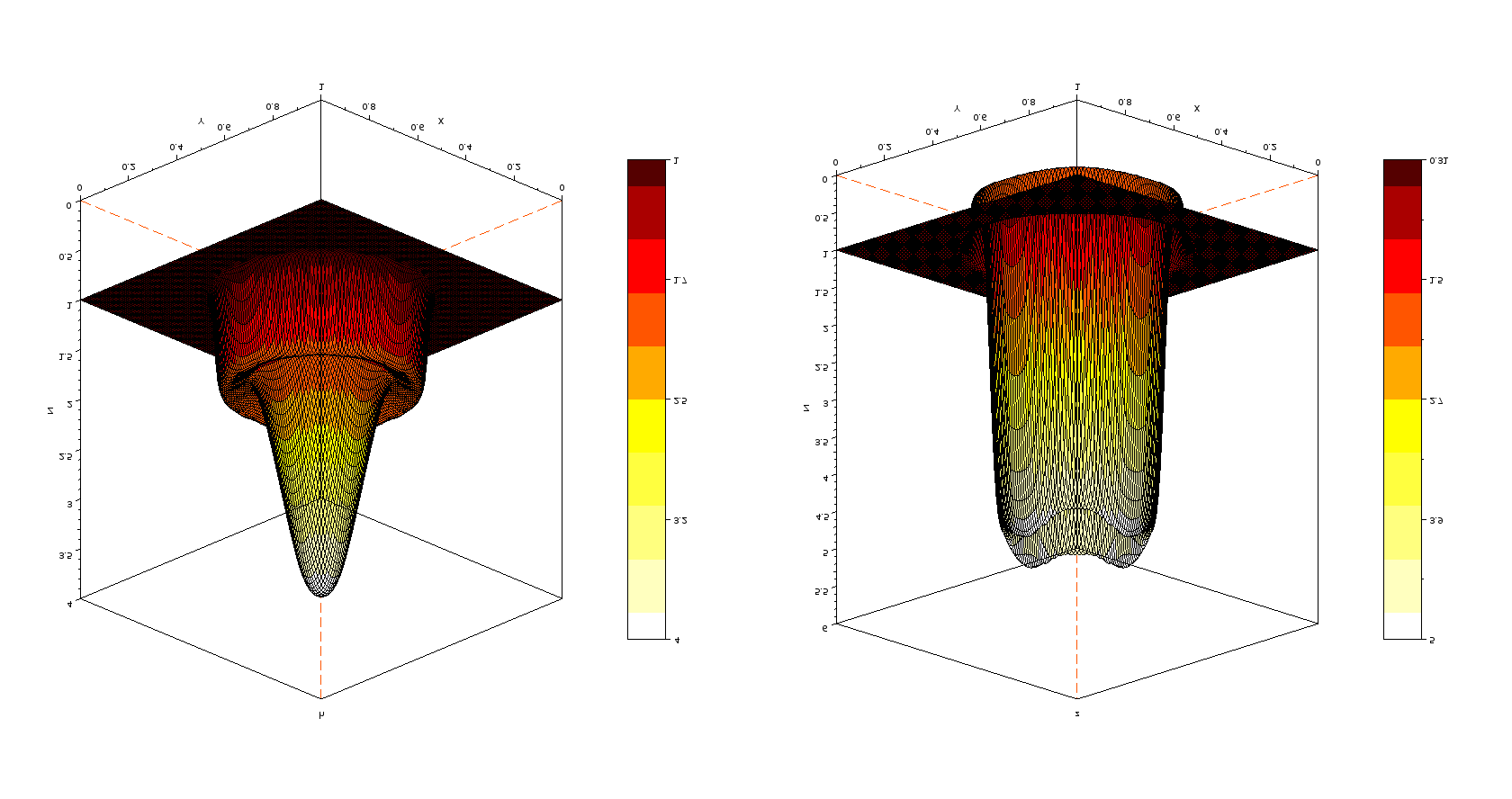}}}
\caption{\label{fig03} Column test case: $H,\Sigma_{zz}$ 
at final time $T=.2$ with a large Weissenberg number $\lambda=1$ 
}
\end{figure}

To emphasize the 2D character of the SVTM and SVUCM models,
we now consider a 
cylindrical version of Stoker test case,
which models the idealized collapse of a fluid column.
A solution for $t\in(0,.2)$ is computed in a square $(x,y)\in[0,1]^2$ starting from the initial condition
$$
(H,U,V,C_{xx},C_{yy},C_{xy},C_{zz}) = 
\begin{cases}
(3,0,0,1,1,0,1) & (x-.5)^2+(y-.5)^2<.2
\\
(1,0,0,1,1,0,1) & (x-.5)^2+(y-.5)^2>.2 
\end{cases}
$$
see Fig.\ref{fig03}.
Note that, a priori, no boundary condition is needed here if we assume our computational domain is a fictitious truncation of the plane $\R^2$
with fictitious boundaries far enough from the initial (circular) discontinuity. 

The main goal of that testcase is usually to see the impact of diverging 
gravity currents on axisymmetric initial conditions.
\mycomment{Self-similar solutions of the axisymmetric shallow-water equations governing converging inviscid gravity currents
can be analytically computed \cite{slim_huppert_2004}}

\begin{figure}[hbtp]
\reflectbox{\rotatebox[origin=c]{180}{\includegraphics[scale=.4,clip,trim = 50 50 50 40]{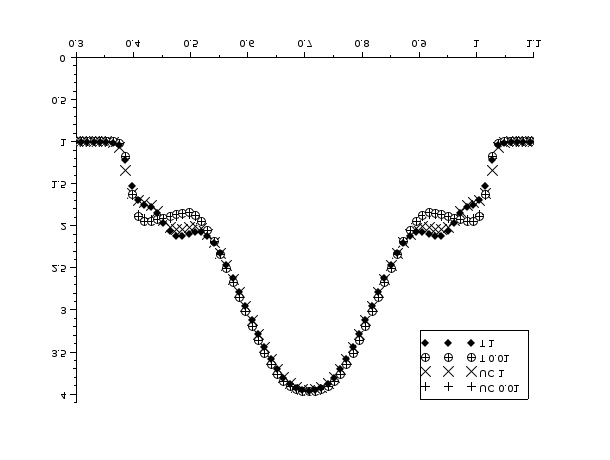}}}
\reflectbox{\rotatebox[origin=c]{180}{\includegraphics[scale=.4,clip,trim = 50 50 50 40]{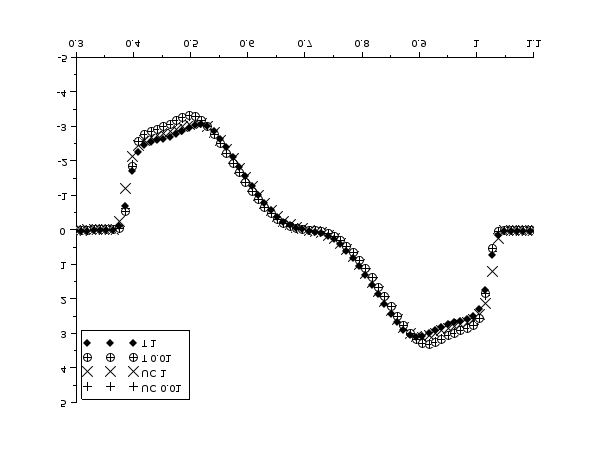}}}
\reflectbox{\rotatebox[origin=c]{180}{\includegraphics[scale=.4,clip,trim = 50 50 50 40]{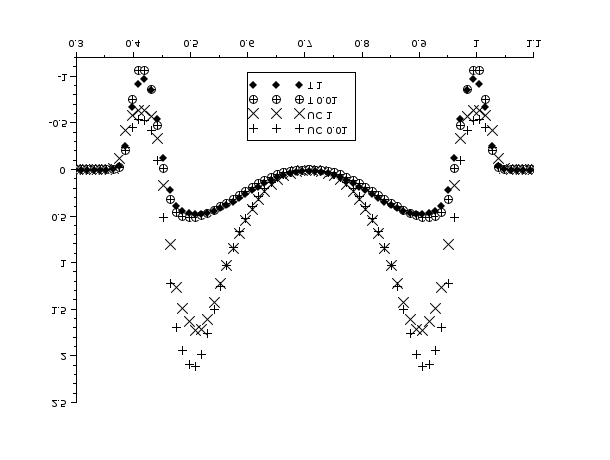}}}
\reflectbox{\rotatebox[origin=c]{180}{\includegraphics[scale=.4,clip,trim = 50 50 50 40]{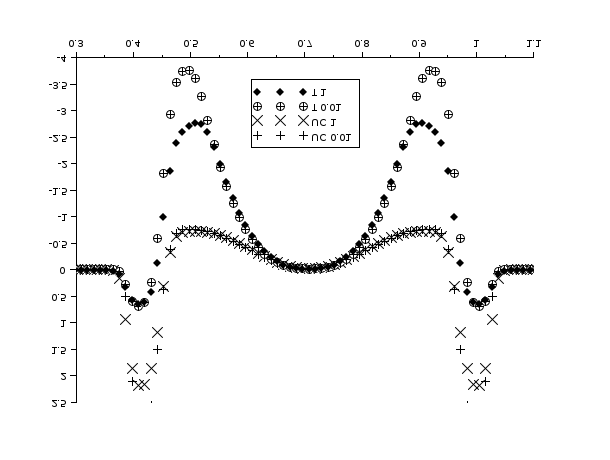}}}
\caption{\label{fig3} Column test case: $H,\bU\cdot\bn,\bSigma\bn\cdot\bn,
\Sigma_{zz}$ 
along diagonal $y=x$ with unit normal $\bn$,
for SVTM and SVUCM with elasticity $G=1$ and $G=.01$, at final time $T=.2$ with a large Weissenberg number $\lambda=1$
(stresses are in $G$ units)}
\end{figure}

Here, we keep the Froude number moderately small $g^{-1/2}=.3$ as before,
and we study the influence of the 
elasticity modulus $G$ for SVTM and SCUVM at a final time $T=.2$.
Note that we can choose the Weissenberg number arbitrarily large,
and the influence of the relaxation-to-equilibrium term is negligible with $\lambda=1\gg T$ for instance.
Comparing the usual hydrodynamical quantities $H,\bU\cdot\bn$ along cross-section $x=y$,
there is little difference between SVTM and SCUVM for the same values $G$.
At large $G=1$, the additional wave clearly shows up in $H$ both for SVTM and SVUCM,
as opposed to the small $G=.001$ case, close to the usual Saint-Venant shallow-water model as expected, see Fig.~\ref{fig3}.

Now, we can also compare $\bC-\bI$ in SVUCM with $\bI-\bC$ in SVTM.
It is then quite striking that the strain deviation from equilibrium (hence the stress) is more important 
for \emph{horizontal} components in SVUCM, 
see the -- most important -- radial component $\bC\bn\cdot\bn-1=\frac1G\bSigma\bn\cdot\bn$ in Fig.~\ref{fig3},
and for \emph{vertical} component $C_{zz}-1=\frac1G\Sigma_{zz}$ in SVTM.
The strain discrepancies in between SVTM and SVUCM should be investigated in the future,
in particular to compare with more standard viscoelastic flow settings with a low Froude number and \emph{shear} forcing
that generates some numerical instabilities in the stress variable (see HWNP below in Section~\ref{sec:ldc}).

\subsection{Lid-driven cavity}
\label{sec:ldc} 

We finally consider a well-known test case for viscous (and viscoelastic) fluid models that involves \emph{stationary} solutions. 
It aims at computing, in a closed square box $(x,y)\in[0,1]^2$ with impermeable walls,
a steady flow satisfying no-slip boundary conditions (i.e. $u=0=v$) at $x=0$, $x=1$, $y=0$ 
and at $y=1$ typically $u=1$, $v=0$ (or a regularized version, see e.g. \cite{Sousa2016129}).
\mycomment{Note that by Kelvin circulation theorem, vorticity is not zero in the cavity, and one vortex at least should develop
which is a priori not a singular structure if the (large-time-limit) velocity field is not purely potential like in the shallow-water standard model.} 

When benchmarking time-dependent (evolutionary) models, stationary solutions are usually computed for large times 
in the hope it becomes close to a limit fixed by the data. 

For viscoelastic fluids however, numerous computer simulations of incompressible creeping flows of Maxwell fluids 
encounter numerical instabilities when elastic stresses increase, see the review about the lid-driven cavity case in \cite{Sousa2016129}.

Precisely, discretizations do not converge to a stationary solution at large Weissenberg number values.
This is one manifestation of the so-called High-Weissenberg-Number-Problem (HWNP).

Various reasons have been invoked to explain the HWNP.
For instance, numerical instabilities may appear 
when the models do not have 
unique solutions anymore.
Now, non-uniqueness may in fact be natural.
Instabilities, i.e. persistent large fluctuations in experimental measures, 
are also sometimes observed physically in settings with (apparently) steady conditions, 
see e.g. the references in \cite{Sousa2016129} as concerns cavity experiments.
But it is not completely clear why the sensitivity of a mathematical model that idealizes the physics
should exactly correspond to the sensitivity of an experimental set-up (unless the model is very good at describing all the physics, 
and its numerical approximation is very accurate).
Indeed, numerical instabilities could also be of a purely mathematical nature, 
and then possibly indicate some imperfection of the (numerical) model at describing the physics, in fact.

That is why, although our aim in the present work is not to ``solve'' the HWNP,
we nevertheless think it is interesting to simulate our new numerical models in more usual conditions for viscoelastic flows
such as the lid-driven cavity, where a HWNP occurs for most existing models.
Indeed, our models somehow enlarge the physical regimes that are usually accessible to (numerical) viscoelastic flows,
with a 
vanishing retardation-time and a non-zero Froude number.

First, 
to obtain conditions that are more usual for viscoelastic flows in a lid-driven cavity,
we choose a low Froude number $g^{-1/2}=10^{-3/2}$ ($g=10^3$), to get closer to the incompressible limit.
Next, we add a viscous component to the stress, 
with a so-called ``solvent viscosity'' $\nu_s=10^{-1}$. 
%
%
Although this does not exactly produce creeping flows 
like in most viscoelastic testcases (our model is evolutionary), the Reynolds number remains quite small $\nu_s^{-1}=10$
so the boundary 
influence is not negligible.

We compute solutions at large times for various values of $G$ and $\lambda$.

%



\begin{figure}[hbtp]
\reflectbox{\rotatebox[origin=c]{180}{\includegraphics[scale=.35,clip,trim = 50 50 50 40]{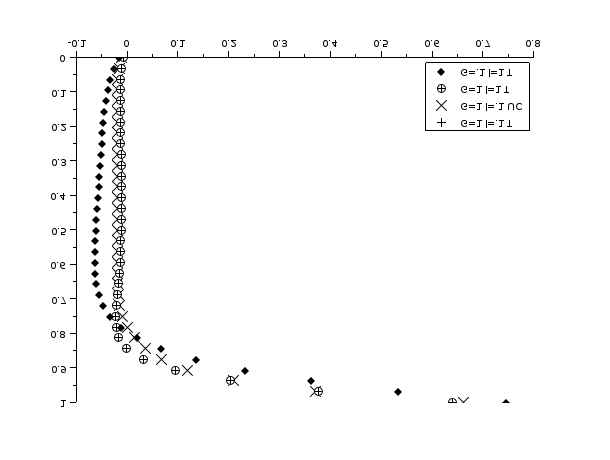}}}
\reflectbox{\rotatebox[origin=c]{180}{\includegraphics[scale=.35,clip,trim = 50 50 50 40]{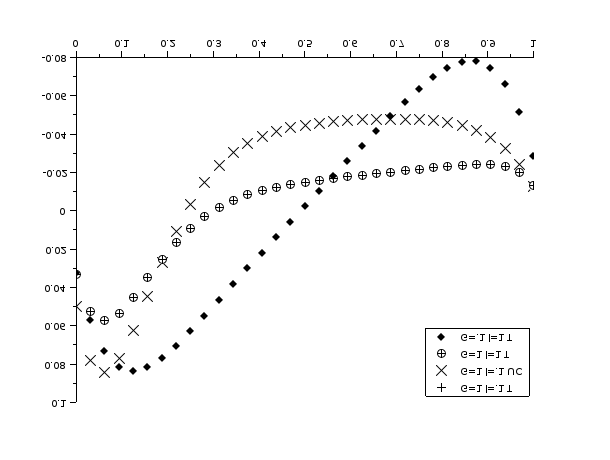}}}
\reflectbox{\rotatebox[origin=c]{180}{\includegraphics[scale=.35,clip,trim = 50 50 50 40]{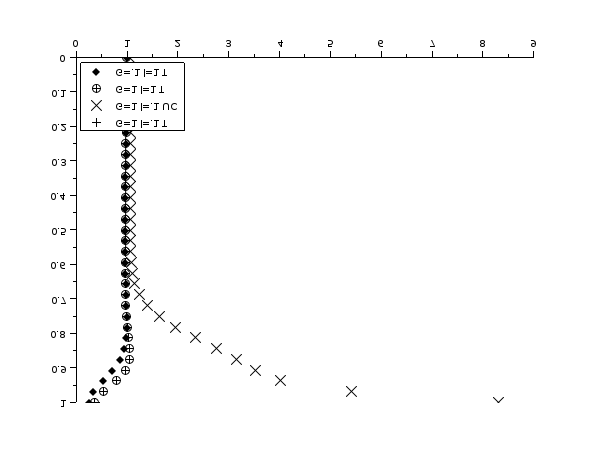}}}
\reflectbox{\rotatebox[origin=c]{180}{\includegraphics[scale=.35,clip,trim = 50 50 50 40]{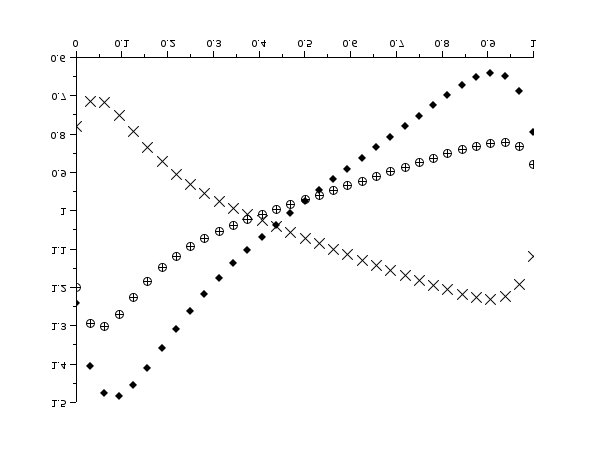}}}
\reflectbox{\rotatebox[origin=c]{180}{\includegraphics[scale=.35,clip,trim = 50 50 50 40]{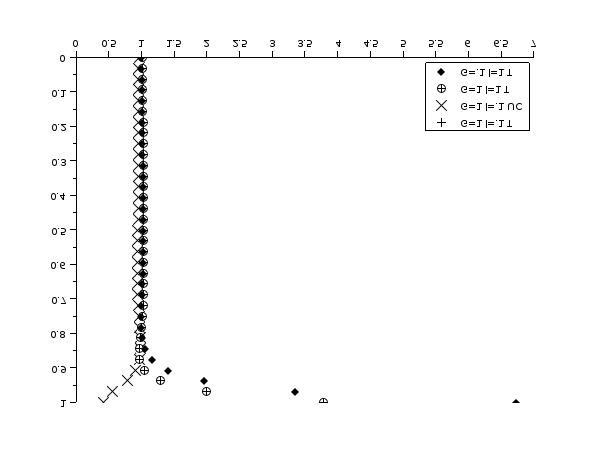}}}
\reflectbox{\rotatebox[origin=c]{180}{\includegraphics[scale=.35,clip,trim = 50 50 50 40]{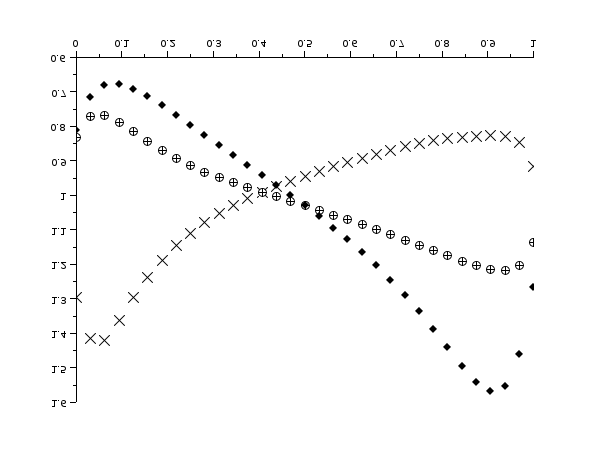}}}
\reflectbox{\rotatebox[origin=c]{180}{\includegraphics[scale=.35,clip,trim = 50 50 50 40]{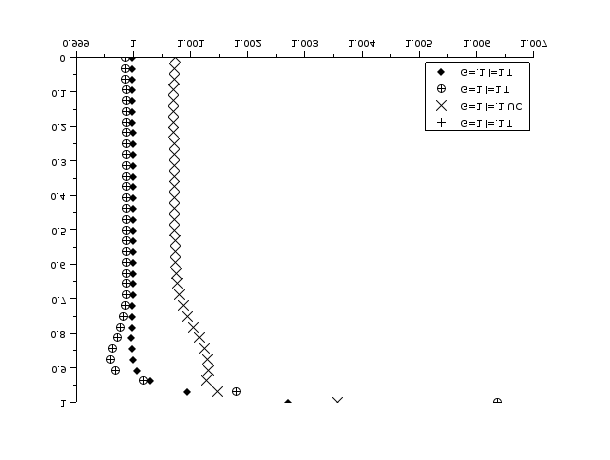}}}
\reflectbox{\rotatebox[origin=c]{180}{\includegraphics[scale=.35,clip,trim = 50 50 50 40]{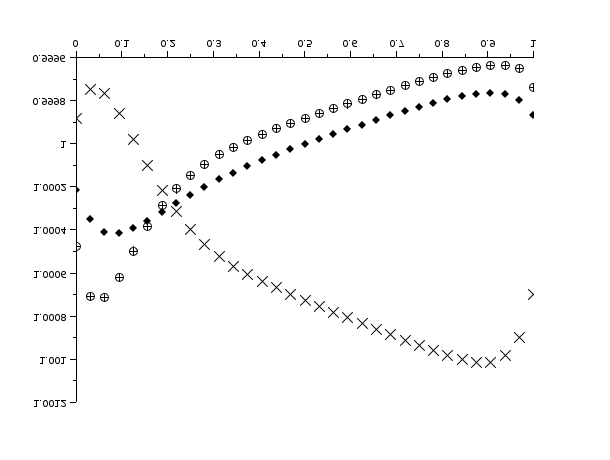}}}
\caption{\label{fig10} Lid-driven cavity at $g=10^3$ with $\nu_s=.1$: cut $x=.5$ (left) and $y=.5$ (right). 
First line: $U$ along $x=.5$ (left) and $V$ along $y=.5$ (right), then $\sx$,$\sy$,$\sz$ (from second to fourth):
for $G=.1,1$, $\lambda=.1,1$, SVTM (T) and SVUCM (UC).}
\end{figure}

In Fig.~\ref{fig10}, we compare at $T=1$ standard quantities of interest ($U$ along $x=.5$, $V$ along $y=.5$\dots) 
computed with a relatively coarse Cartesian mesh of $33\times 33=1089$ cells.
First, whereas $H$, $U$ and $V$ are hardly different for SVTM and SVUCM,
the viscoelastic stress are different: $G(\bI-\bC)$ in SVUCM does not exactly match $G(\bC-\bI)$ in SVTM,
and the components of the two stress tensors can be quite larger for SVUCM than SVTM
(although they both have similar variations around zero), see $C_{xx}$ along $x=.5$ in Fig.~\ref{fig10}.
Second, the stationary solution seems determined by $G$ only,
and not by our Weissenberg number $\lambda$ unlike the usual ``creeping flow'' solutions
typically computed at fixed $\beta:=\nu_s/(\nu_s+\nu_p)=.5$ with a ``polymer viscosity'' $\nu_p:=G\lambda$.

\begin{figure}[hbtp]
\reflectbox{\rotatebox[origin=c]{180}{\includegraphics[scale=.35,clip,trim = 50 50 50 40]{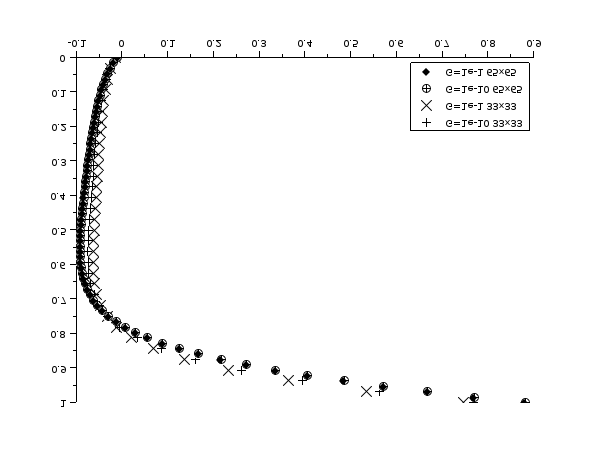}}}
\reflectbox{\rotatebox[origin=c]{180}{\includegraphics[scale=.35,clip,trim = 50 50 50 40]{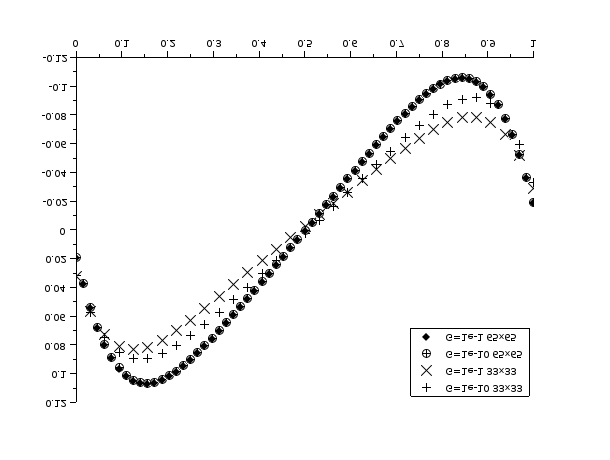}}}
\reflectbox{\rotatebox[origin=c]{180}{\includegraphics[scale=.35,clip,trim = 50 50 50 40]{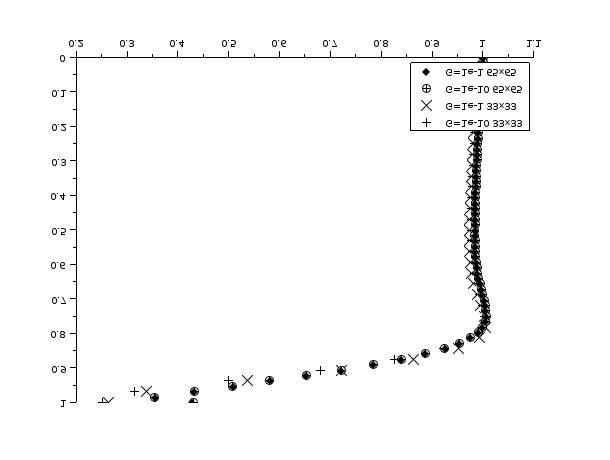}}}
\reflectbox{\rotatebox[origin=c]{180}{\includegraphics[scale=.35,clip,trim = 50 50 50 40]{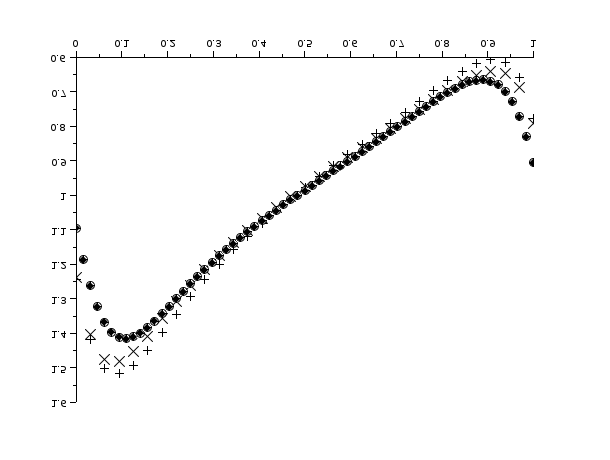}}}
\reflectbox{\rotatebox[origin=c]{180}{\includegraphics[scale=.35,clip,trim = 50 50 50 40]{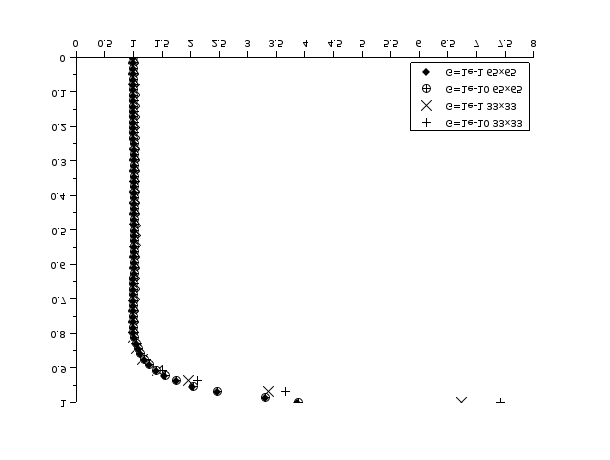}}}
\reflectbox{\rotatebox[origin=c]{180}{\includegraphics[scale=.35,clip,trim = 50 50 50 40]{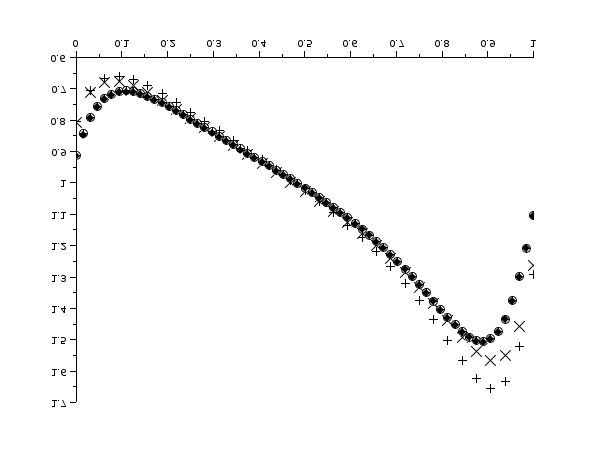}}}
\reflectbox{\rotatebox[origin=c]{180}{\includegraphics[scale=.35,clip,trim = 50 50 50 40]{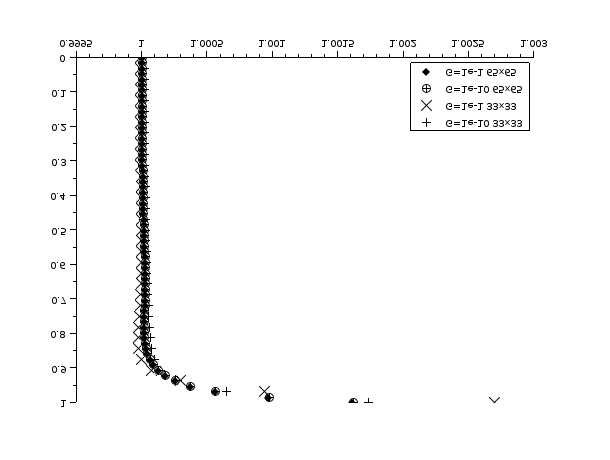}}}
\reflectbox{\rotatebox[origin=c]{180}{\includegraphics[scale=.35,clip,trim = 50 50 50 40]{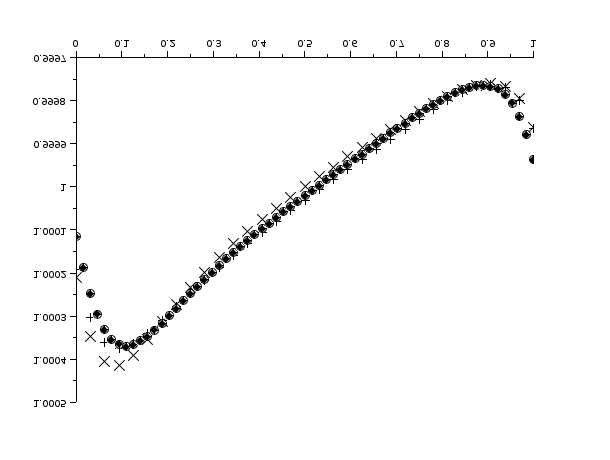}}}
\caption{\label{fig12} Lid-driven cavity at $g=10^3$ with $\nu_s=.1$: cut $x=.5$ (left) and $y=.5$ (right). 
First line: $U$ along $x=.5$ (left) and $V$ along $y=.5$ (right), then $\sx$,$\sy$,$\sz$ (from second to fourth):
for $\lambda=1$, SVTM, $G=.1,1^{-10}$ and two meshes.}
\end{figure}

\begin{figure}[hbtp]
\reflectbox{\rotatebox[origin=c]{180}{\includegraphics[scale=.35,clip,trim = 50 40 50 40]{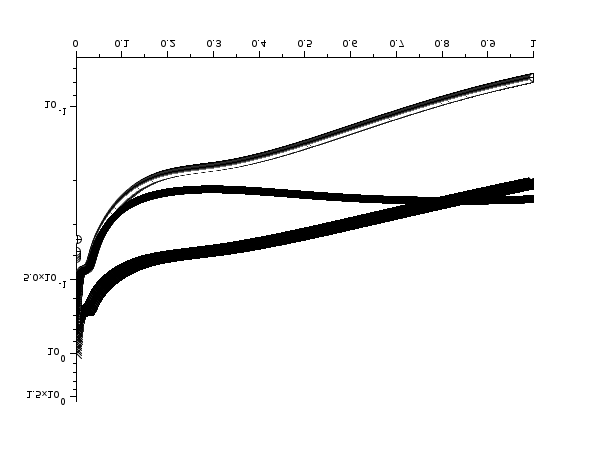}}}
\reflectbox{\rotatebox[origin=c]{180}{\includegraphics[scale=.35,clip,trim = 50 40 50 40]{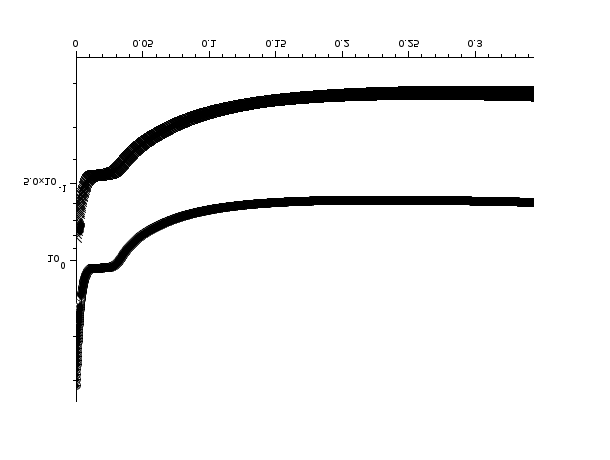}}}
\caption{\label{fig11} Lid-driven cavity at $g=10^3$ with $\nu_s=.1$. Left:
$\ell^1$ norms of the differences between two successive solutions as a function of time (our stationarity criterium)
for the four cases of Fig.~\ref{fig10} with $33\times 33=1089$ cells,
namely $G=1$,$\lambda=.1$ UC then $G=.1$,$\lambda=1$ T and $G=1$,$\lambda=.1,1$ T (superimposed) from top to bottom. 
Right: stagnation reaches a plateau higher with $65\times 65=4225$ than with $33\times 33=1089$ cells,
for $G = 10^{-1}$ and $G=10^{-10}$ (superimposed). 
}
\end{figure}

Note however that the time-dependent numerical solutions of Fig.~\ref{fig10} are quickly stationary in time for the smallest value $G=.1$,
and seem limited in accuracy due to the presence of viscosity.  
In Fig.~\ref{fig11}, we refine the mesh to $65\times 65=4225$ cells:
the $\ell^1$ norm of the differences between two successive solutions as a function of time (our stationarity criterium)
stagnates at a higher level, 
a plateau which is the same for all values of $G$ smaller than $.1$
(see Fig.~\ref{fig12} for the converged solutions when $G = 10^{-1}$ and $G=10^{-10}$).
So our numerical experiment seems actually 
interesting only for large enough values of $G$, 
when shearing becomes more difficult
and when strain is not so large  
(see the variations of $U$ and $V$ with $G$ in Fig.~\ref{fig10}). %

\begin{figure}[hbtp]
\reflectbox{\rotatebox[origin=c]{180}{\includegraphics[scale=.35,clip,trim = 50 40 50 40]{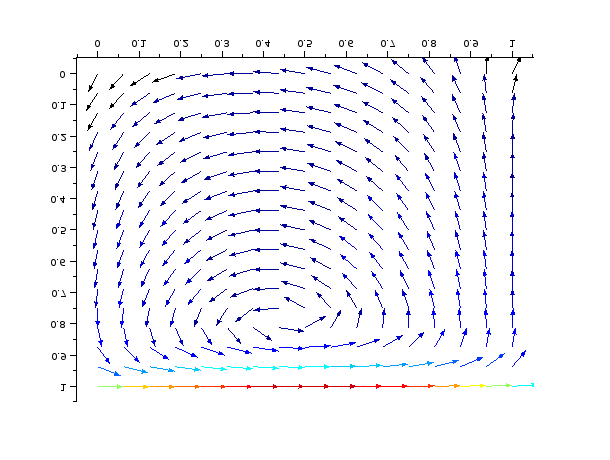}}}
\reflectbox{\rotatebox[origin=c]{180}{\includegraphics[scale=.35,clip,trim = 50 40 50 40]{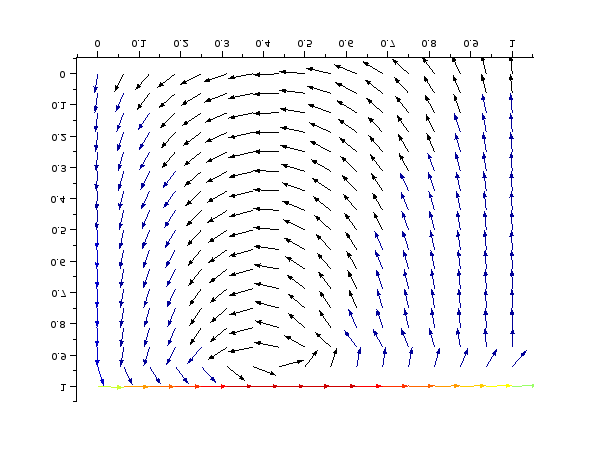}}}
\reflectbox{\rotatebox[origin=c]{180}{\includegraphics[scale=.35,clip,trim = 50 40 50 40]{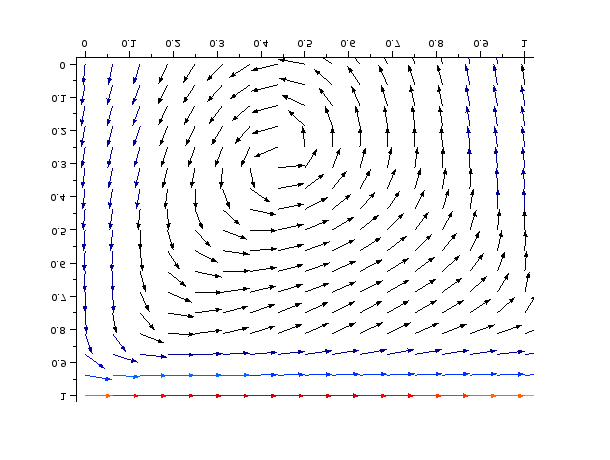}}}
\reflectbox{\rotatebox[origin=c]{180}{\includegraphics[scale=.35,clip,trim = 50 40 50 40]{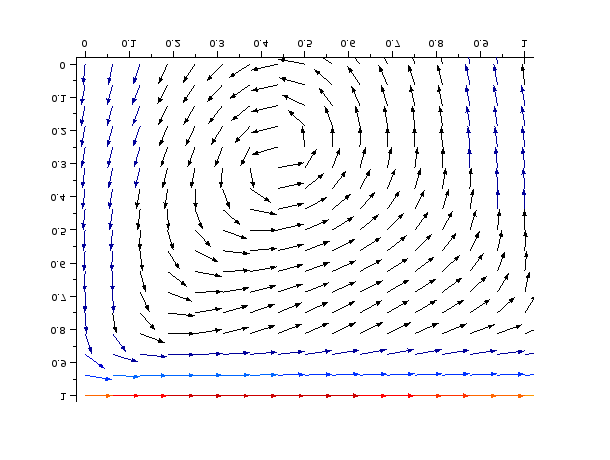}}}
\caption{\label{fig13} Lid-driven cavity at $g=10^3$: SVTM velocity vector fields $\bU=(U,V)$ at $T=1$
for $\nu_s=.1,.01$ (left/right) and $G=1,10$ (top/bottom).
}
\end{figure}

\begin{figure}[hbtp]
\reflectbox{\rotatebox[origin=c]{180}{\includegraphics[scale=.35,clip,trim = 50 40 50 40]{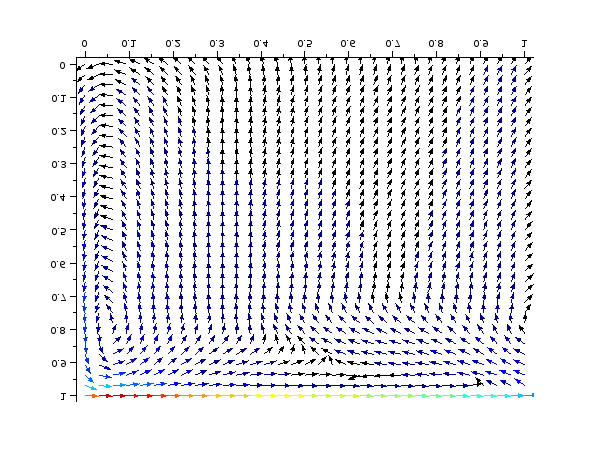}}}
\reflectbox{\rotatebox[origin=c]{180}{\includegraphics[scale=.35,clip,trim = 50 40 50 40]{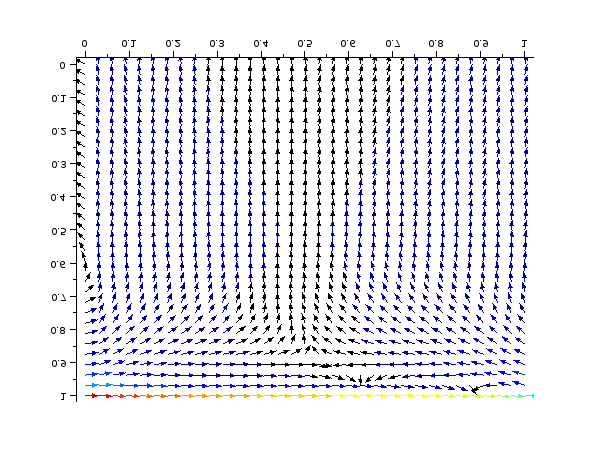}}}
\reflectbox{\rotatebox[origin=c]{180}{\includegraphics[scale=.35,clip,trim = 50 40 50 40]{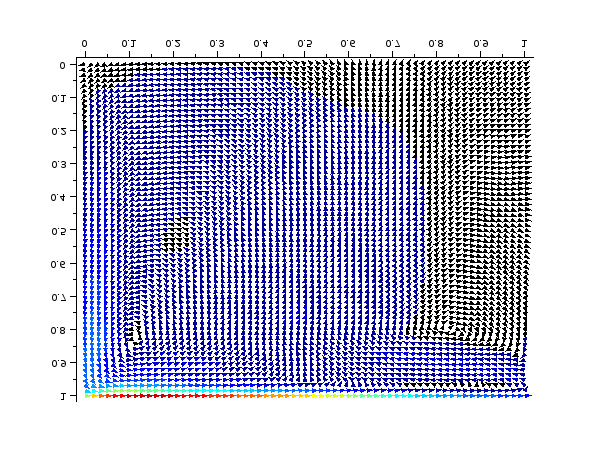}}}
\reflectbox{\rotatebox[origin=c]{180}{\includegraphics[scale=.35,clip,trim = 50 40 50 40]{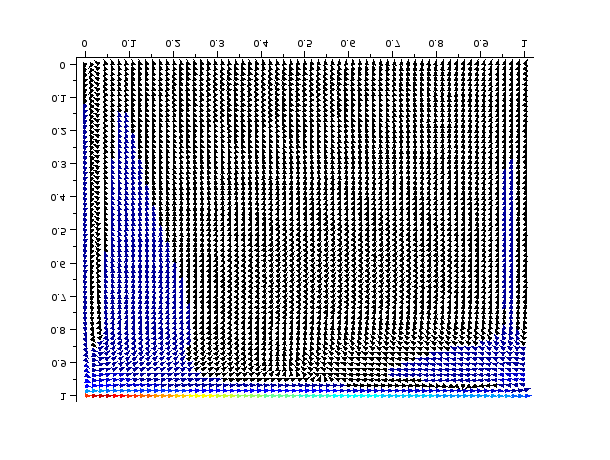}}}
\caption{\label{fig14} Lid-driven cavity at $g=10^3$: SVUCM velocity vector field $\bU=(U,V)$
at $33\times 33=1089$ and $65\times 65=4225$ points (top/bottom) for $G=1,10$ (left/right) and $\nu_s=10^{-2}$, $\lambda=.1$.
}
\end{figure}

Now, for larger $G$, we do observe convergence in time to stationary states without reaching a plateau both for SVTM and SVUCM.
However, SVTM and SVUCM solutions now strongly differ.
Solutions to SVT and SVUCM almost coincide at $\nu_s=.1$ and $G=1$ and this is easily seen from the velocity vector fields $\bU=(U,V)$
(see e.g. Fig.~\ref{fig13}).
Then, if we increase $G$ to 10, large-time SVTM simulations seem to converge (in time and space) to solutions
with only one main vortex, which is only slightly deformed and influenced by $\nu_s,\lambda$ 
(see in Fig.~\ref{fig13}).
On the contrary, SVUCM converge to a different type of solution which is also captured when $\nu_s$ is smaller, see Fig.~\ref{fig14}.

\section{Conclusion}

In this work, motivated by the need for better numerical models of viscoelastic flows,
we first have derived new hyperbolic models in the framework of shallow free-surface flows proposed by Saint-Venant.
This extends Saint-Venant 2D shallow-water model to Maxwell fluid and is a continuation of our 1D work \cite{bouchut-boyaval-2013}.

One model coincides with the zero-retardation Oldroyd-B case which had been obtained somewhat differently in \cite{bouchut-boyaval-2015}
(without a precise study of solutions like here).
The other suggests one to use a viscoelastic equation for a conformation tensor with a time-rate different than what is usually done in the literature
(i.e. the Gordon-Schowalter derivatives).

Next, we have also proposed Finite-Volume (FV) discretizations that preserve the essential properties of the new models:
mass and momentum conservation, plus a free-energy dissipation.
Numerical simulations have been performed with the FV schemes that phenomenologically prove the physical soundness of the model
in simple free-shear flows.

Quantitative evaluations in the lid-driven cavity testcase 
also show the interest of the approach to 
investigate more realistic strongly-sheared viscoelastic flows.
Our numerical scheme should however still be improved to better understand standard benchmarks with large strains (and then with a HWNP, quite often),
where differences between the two models (SVTM and SVUCM) with two different time-rates may be important.
In particular, accuracy should be improved in the low Froude regime.
This will be the object of future work.


\bibliographystyle{amsplain}

\providecommand{\bysame}{\leavevmode\hbox to3em{\hrulefill}\thinspace}
\providecommand{\MR}{\relax\ifhmode\unskip\space\fi MR }
\providecommand{\MRhref}[2]{%
  \href{http://www.ams.org/mathscinet-getitem?mr=#1}{#2}
}
\providecommand{\href}[2]{#2}

\end{document}